\documentclass[onefignum,onetabnum]{siamart190516}

\usepackage{graphicx}

\usepackage{geometry}
\newgeometry{
	textwidth=13.3cm,
        textheight=23cm,
      }
\usepackage{amsmath}
\usepackage{amssymb}
\usepackage{algorithmic}

\usepackage{enumerate}
\usepackage{url}

\newcommand{\FB}{T_{\rm{FB}}}

\newcommand{\id}{\mathrm{Id}}

\newcommand{\proj}{\Pi}
\newcommand{\reals}{\mathbb{R}}

\newcommand{\hilbert}{\mathcal{H}}

\newcommand{\weakto}{\rightharpoonup}

\newcommand{\fix}{\mathrm{fix}}
\newcommand{\zero}{\mathrm{zer}}
\newcommand{\dist}{\mathrm{dist}}

\newcommand{\Fejer}{Fej{\'e}r}
\newcommand{\gph}{\mathrm{gph}}
\newcommand{\naturals}{\mathbb{N}}
\newcommand{\interior}{{\rm{int}}}

\newcommand{\dom}{{\rm{dom}}}

\newcommand{\ran}{{\rm{ran}}}

\newcommand{\PD}{\mathcal{P}(\hilbert)}
\newcommand{\neighborU}{\mathcal{U}}
\newcommand{\neighborV}{\mathcal{V}}
\newcommand{\ordth}{'th~}

\DeclareMathOperator*{\argmin}{argmin}
\DeclareMathOperator*{\minimize}{minimize}

 \newtheorem{prop}{Proposition}
\newtheorem{lem}{Lemma}
 \newtheorem{cor}{Corollary}
 \newtheorem{ass}{Assumption}
 \newtheorem{thm}{Theorem}
 \newtheorem{defin}{Definition}
 \newtheorem{rem}{Remark}
 \newenvironment{pf}{\smallbreak\noindent{\it Proof. }}{\hfill$\Box$\smallbreak}

\xdef\one{1}
\xdef\showalgs{1} 
\xdef\showps{1}  

\newcommand{\takeout}[1]{}

\title{Nonlinear Forward-Backward Splitting\\with Projection Correction\thanks{Financial support from the Swedish Research Council is gratefully acknowledged.}}
\author{Pontus Giselsson}
\date{}

\begin{document}
\maketitle

\begin{abstract}
  We propose and analyze a versatile and general algorithm called nonlinear forward-backward splitting (NOFOB). The algorithm consists of two steps; first an evaluation of a {\emph{nonlinear}} forward-backward map followed by a relaxed projection onto the separating hyperplane it constructs. The key of the method is the nonlinearity in the forward-backward step, where the backward part is based on a {\emph{nonlinear}} resolvent construction that allows for the kernel in the resolvent to be a nonlinear single-valued maximal monotone operator. This generalizes the standard resolvent as well as the Bregman resolvent, whose resolvent kernels are gradients of convex functions. This construction opens up for a new understanding of many existing operator splitting methods and paves the way for devising new algorithms. In particular, we present a four-operator splitting method as a special case of NOFOB that relies nonlinearity and nonsymmetry in the forward-backward kernel. We show that forward-backward-forward splitting (FBF), forward-backward-half-forward splitting (FBHF), asymmetric forward-backward-adjoint splitting (AFBA) with its many special cases, as well as synchronous projective splitting are special cases of the four-operator splitting method and hence of NOFOB. We also show that standard formulations of FB(H)F use smaller relaxations in the projections than allowed in NOFOB. Besides proving convergence for NOFOB, we show linear convergence under a metric subregularity assumption, which in a unified manner shows (in some cases new) linear convergence results for its special cases. 
  
\end{abstract}

\begin{keywords}
  Monotone inclusions, Nonlinear resolvent, Forward-backward splitting, Forward-backward-forward splitting, Four-operator splitting
\end{keywords}
  
\begin{AMS}
 90C25, 65K05, 90C30
\end{AMS}

\section{Introduction}

We consider maximal monotone inclusion problems of the form
\begin{align}
  0\in Ax+Cx,
  \label{eq:NOFOBprob}
\end{align}
where $A:\hilbert\to 2^\hilbert$ is maximally monotone, $C:\hilbert\to\hilbert$ is cocoercive, and $\hilbert$ is a real Hilbert space. This problem is ubiquitous in engineering fields as it comprises problems from optimization, variational analysis, and game theory. 
We present a flexible algorithm called nonlinear forward-backward splitting (NOFOB) for solving \eqref{eq:NOFOBprob}. It has many algorithms, including the classical forward-backward splitting (FBS) \cite{Mercier1979Lectures,Passty1979Ergodic,Combettes2005Signal}, as special cases. The first step in NOFOB is a nonlinear forward-backward step. This is a novel construction that generates a separating hyperplane between the current point and the solution set. A relaxed projection onto this separating hyperplane finishes one iteration.

The nonlinear forward-backward map applied to $A$ and $C$ is defined as
\begin{align*}
\FB:=(M+A)^{-1}\circ (M-C),
\end{align*}
where $M:\hilbert\to\hilbert$ is a (possibly) nonlinear strongly monotone and Lipschitz continuous operator. The nonlinearity of the forward-backward map is in the kernel $M$ of the backward part $(M+A)^{-1}\circ M$, that we call a {\emph{nonlinear resolvent}}.
If $M=\id$, the standard forward-backward map is recovered. If instead $M=\nabla g$ for a differentiable convex function $g$, the Bregman forward-backward map is recovered. The Bregman resolvent has been thoroughly studied in the literature, see, e.g., \cite{Eckstein1993Nonlinear,Bauschke2003Bregman}, and later also with forward steps \cite{Bauschke2017ADescent,Lu2018Relatively}. That we allow for arbitrary maximal monotone kernels $M$ is the key that, e.g., allows us to cast FBF and FBHF as special cases of our method.

We utilize the versatility that comes from the possibility to use nonlinear and nonsymmetric kernels to present a novel four-operator splitting method as a special case of NOFOB. It solves problems of the form
\begin{align*}
  0\in Bx+Dx+Ex+Kx,
\end{align*}
where $B:\hilbert\to 2^\hilbert$, $D:\hilbert\to\hilbert$ is Lipschitz continuous, $B+D$ is maximally monotone, $E:\hilbert\to\hilbert$ is cocoercive, and $K$ is linear skew-adjoint, which is an important property in primal-dual methods. The method is obtained from NOFOB by letting $A=B+D+K$, $C=E$, and $M=Q-D-K$, where $Q$ is a strongly monotone and Lipschitz continuous operator. This choice of $M$ effectively moves the single-valued operators $D$ and $K$ from the backward part to the forward part, since
\begin{align*}
  \hat{x}&=(M+A)^{-1}(M-C)x \\&=((Q-D-K)+(B+D+K))^{-1}(Q-D-K-E)x \\
                   &=(Q+B)^{-1}(Q-(D+K+E))x.
\end{align*}
This results in that $D$ and $K$ need only be accessed via forward evaluation.
As in NOFOB, it is not enough to iterate this nonlinear forward-backward map as the resulting sequence may diverge. The four-operator splitting method uses the same projection correction step as in NOFOB that performs a relaxed projection onto the separating hyperplane that is generated by the nonlinear forward-backward step. The forward-backward step evaluates $M$ at $x$ and the projection requires an additional application of $M$ at the point $\hat{x}$. Therefore, $D$ and $K$ need to be evaluate twice per iteration, while $E$ only needs to be evaluated once. 

We will show that the four-operator splitting method, and hence NOFOB, has many special cases in the literature. By letting $K=E=0$, $Q=\gamma^{-1}\id$, where $\gamma>0$ is a step-size, and using a conservative relaxation factor in the relaxed projection, we arrive at a standard formulation of forward-backward-forward splitting \cite{Tseng2000AModified} (which is equivalent to the extra-gradient method (EG) \cite{Korpelevich1976TheExtragradient} if also $B=0$). FBF and EG have recently gained momentum in the literature, partly due to their stabilizing properties when training generative adversarial networks, GANs \cite{Goodfellow2014Generative,Gidel2019AVariational}. Forward-backward-half-forward splitting~\cite{BricenoArias2018Forward} is obtained from the four-operator splitting method by letting $K=0$, $Q=\gamma^{-1}\id$, and using a conservative relaxation factor in the relaxed projection. The asymmetric forward-backward-adjoint splitting method (AFBA) in \cite{Latafat2017Asymmetric,Latafat2018PrimalDual}, is obtained by letting $D=0$ and $Q=P+G$, where $P$ is self-adjoint positive definite and $G$ is linear skew-adjoint. Therefore, the lists of special cases that are cataloged in \cite{Latafat2017Asymmetric,Latafat2018PrimalDual}, such as Solodov and Tseng~\cite{Solodov1996Modified}. Douglas-Rachford~\cite{Lions1979Splitting}, ADMM with dual step-length 1~\cite{Gabay1976ADual,Glowinski1975Sur}, and \cite{Drori2015ASimple,BricenoArias2011AMonotone,He2012Convergence,Latafat2019ANew}, are also special cases of this algorithm. We will also show that a synchronous version of projective splitting \cite{Eckstein2008AFamily,Combettes2018Asynchronous} (that has been further analyzed in \cite{Johnstone2018Projective,Johnstone2019Convergence}) is a special case with $D=E=0$ and a positive definite $Q$. Since $K$ and $D$ are treated similarly in the four-operator splitting method, this reveals that projective splitting has the same algorithmic structure as FBF.

The standard forward-backward splitting method (along with its special cases, e.g., Chambolle-Pock~\cite{Chambolle2011AFirstOrder} and V{\~u}-Condat~\cite{Condat2013APrimalDual,Vu2013ASplitting}) is a special case of NOFOB and the four-operator splitting method. If $M=\gamma^{-1}\id$ with step-size $\gamma>0$, we get the standard forward-backward step $\FB=(\id+\gamma A)^{-1}(\id-\gamma C)$. In this case, we can select projection metric so that the projection point is the same as the result of the forward-backward step. Therefore, the projection step becomes redundant, the second application of $M$ is not needed, and NOFOB reduces to FBS. This implies that NOFOB automatically adapts to nonlinearity or nonsymmetry of $M$. If $M$ is nonlinear or nonsymmetric, NOFOB makes use of the projection correction step, that needs a second evaluation of $M$, to guarantee convergence, while if $M$ is linear, self-adjoint, and positive definite the correction step is redundant and the second application of $M$ is avoided.

We present and prove convergence of NOFOB with iteration dependent kernel $M$. We first show that the introduced nonlinear forward-backward map indeed creates a separating hyperplane between the current iterate and the solution set. The separating hyperplane defines a halfspace $H_k$ that contains the solution set but not the current point. We call this halfspace a separating halfspace. The relaxed projection onto $H_k$ gives \Fejer~monotonicity w.r.t. to the solution set, i.e., the distance to the solution set is nonincreasing between iterations. The second part shows that the cuts are deep enough, i.e., that the separating halfspace is far enough from the current point for the algorithm to not stall before reaching the solution set. We also show linear convergence of NOFOB in finite-dimensional settings under a metric subregularity assumption. These results provide a unified convergence analysis for all the special cases and in particular linear convergence results for forward-backward-half-forward splitting \cite{BricenoArias2018Forward} and synchronous projective splitting \cite{Eckstein2008AFamily,Combettes2018Asynchronous} that appear to be new.

\subsection{Related Work}

The nonlinear forward-backward algorithm is related to, and generalizes, many methods in the literature. The Hybrid Projection-Proximal Point Method (HPPPM) and variations are proposed and analyzed in a sequence of papers \cite{Solodov1999AHybrida,Solodov1999AHybridb,Solodov2001AUnified}. The algorithms are variations of Rockafellar's proximal point algorithm \cite{Rockafellar1976Monotone} that allow for specific inexact resolvent updates. These updates are followed by a correction step in the form of a projection onto a separating hyperplane, see \cite{Solodov2001AUnified} that provides a unified treatment of \cite{Solodov1999AHybrida,Solodov1999AHybridb}. The kernel $M$ used in their resolvent step is always the identity. This has been generalized to Bregman operators in \cite{Solodov2000AnInexact}, i.e., where $M$ is the gradient of a differentiable convex function. This results in an approximate Bregman resolvent method. Such methods have also been analyzed in \cite{Eckstein1998Approximate} without the correction step. Our method is different from the above in that we allow for arbitrary nonlinear strongly monotone kernels $M$ and that we have an additional cocoercive term in the model. It was shown in \cite{Solodov1999AHybridb} that FBF is a special case of the {\emph{inexact}} proximal point method with correction proposed in \cite{Solodov1999AHybridb}. By allowing for a nonlinear $M$ in the resolvent in NOFOB, we can cast FBF as an {\emph{exact}} special case.

Similar algorithms based on the separate and project principle have been proposed for equilibrium problems and variational inequalities (that are special instances of monotone inclusion problems, see \cite{Combettes2005Equilibrium,Aoyama2008Maximal}) in \cite{Konnov1993Combined,Konnov1999ACombined,Konnov2005Combined,Konnov2006Combined}. The separate and project principle is in these works referred to as a {\emph{combined relaxation approach}}. The methods call a separation oracle that provides a sufficiently deep separation followed by a relaxed projection onto the separating hyperplane. Many specific instances of separation oracles based on resolvent evaluations have been suggested for variational inequalities and equilibrium problems in \cite{Konnov1993Combined,Konnov1999ACombined,Konnov2005Combined,Konnov2006Combined}. NOFOB falls within the same very general framework of separate and project methods, but generalizes the above method instances in that it allows for arbitrary maximal monotone inclusion problems, infinite dimensional Hilbert spaces, as well as a cocoercive term in the model. 

Less similar, but still related methods based on projection onto separating hyperplanes have been proposed for finding a point in the intersection of a finite number of convex sets \cite{Kiwiel1995BlockIterative} and for finding a common fixed-point of firmly nonexpansive mappings \cite{Kiwiel1997Surrogate}. These methods are generalized in \cite{Combettes2000Strong} that considers minimization of a uniformly convex function over a finite number of convex sets. 


The introduced nonlinear forward-backward map reduces to a nonlinear resolvent map in absence of cocoercive term $C$. Properties of the resolvent map have been investigated in \cite{Bauschke2010General}, but with no mentioning of algorithms. This has been further investigated in \cite{Bui2019Warped}, which has been developed in parallel with this manuscript. The two works (\cite{Bui2019Warped} and this manuscript) are truly independent as the first respective preprints were uploaded during the same arXiv upload slot. The papers share the resolvent construction (with different names), but propose and analyze different general algorithmic schemes based on the resolvent, and consider, to a large extent, different operator splitting methods to be special cases of the respective schemes. Therefore, the work in \cite{Bui2019Warped} is a nice complement to our work.

Finally, we note a relation to the four-operator splitting method in \cite{MarquesAlves2018Iteration}, whose problem structure also consists of four operators, but differs in that the operators have other properties.

\section{Preliminaries}
\label{sec:Prelim}

In this section, we collect notation, definitions, and simple auxiliary results. We let $\reals$ be the set of real numbers and $\reals_{+}$ denote the set of positive real numbers. Further, $\hilbert$ and $\mathcal{G}_i$ denote real Hilbert spaces. Their inner products and induced norms are denoted $\langle\cdot,\cdot\rangle$ and $\|\cdot\|$ respectively. We let $\PD$ be the set of bounded linear self-adjoint positive definite operators on $\hilbert$. For $P\in\PD$, we define $\|x\|_P=\sqrt{\langle x,Px\rangle}$.
The operator norm of bounded linear operator $L$ induced by $\|\cdot\|_P$ is $\|L\|_P:=\sup_{x\in\hilbert:\|x\|_P\leq 1}\|Lx\|_{P}$.
We denote by $\lambda_{\min}(P),\lambda_{\max}(P)\in\reals_+$ the numbers that satisfy $\lambda_{\min}(P)\|x\|^2\leq\|x\|_P^2\leq\lambda_{\max}(P)\|x\|^2$ for all $x\in\hilbert$. The $\|\cdot\|_S$-norm distance from $x\in\hilbert$ to a nonempty closed convex set $C\subset\hilbert$ is defined as $\dist_S(x,C):=\inf_{y\in C}\|x-y\|_S$ and satisfies $\dist_S(x,C)=\|x-\Pi_C^Sx\|_S$, where $\Pi_C^Sx$ is the projection $\Pi_C^Sx=\argmin_{y\in C}\|x-y\|_S$.
Further, we use the convention that $\tfrac{0}{0}=0$ and $\tfrac{\alpha}{0}=\infty$ for all $\alpha\in\reals_+$.

The powerset of $\hilbert$ is denoted $2^\hilbert$.
The graph of $A:\hilbert\to 2^\hilbert$ is $\gph(A)=\{(x,u):u\in Ax\}$.
  An operator $A:\hilbert\to 2^\hilbert$ is {\emph{monotone}} if
    $\langle u-v,x-y\rangle\geq 0$
  for all $(x,u),(y,v)\in\gph(A)$. It is {\emph{$\sigma$-strongly monotone}} w.r.t. $\|\cdot\|_P$ with $P\in\PD$ if $\sigma>0$ and
    $\langle u-v,x-y\rangle\geq\sigma\|x-y\|_{P}^2$
  for all $(x,u),(y,v)\in\gph(A)$.
  An operator $T:\hilbert\to\hilbert$ is {\emph{$\beta$-cocoercive}} w.r.t. $\|\cdot\|_S$ with $S\in\PD$ if $\beta\geq 0$ and
    $\langle Tx-Ty,x-y\rangle\geq\beta\|Tx-Ty\|_{S^{-1}}^2$
  for all $x,y\in\hilbert$, and {\emph{$L$-Lipschitz continuous}} w.r.t. $\|\cdot\|_S$ with $L\geq 0$ if
    $\|Tx-Ty\|_{S^{-1}}\leq L\|x-y\|_S$
  for all $x,y\in\hilbert$. 
  Finally, an operator $K:\hilbert\to\hilbert$ is {\emph{skew}} if
    $\langle Kx-Ky,x-y\rangle=0$
  for all $x,y\in\hilbert$. If $K$ in addition is bounded and linear, it satisfies $-K=K^*$. We call such operators {\emph{linear skew-adjoint}}.

We finish by stating the following simple results for ease of reference.
\takeout{\begin{prop}
  An operator $T$ is $\beta$-cocoercive w.r.t. $\|\cdot\|$ with $\beta>0$ if and only if $T=\tfrac{1}{2\beta}(\id+N)$ for some nonexpansive (w.r.t. $\|\cdot\|$) operator $N$.
  \label{prop:coco_nonexp}
\end{prop}
\begin{pf}
$T$ is $\beta$-cocoercive if and only if $\beta T$ is firmly nonexpansive \cite[Definition~4.10]{Bauschke2017Convex}, i.e., $\beta T=\tfrac{1}{2}(\id+N)$ for some nonexpansive $N$ \cite[Proposition~4.4]{Bauschke2017Convex}. Dividing by $\beta$ gives the result.
\end{pf}}
\takeout{\begin{prop}
  An operator $T:\hilbert\to\hilbert$ that is $\tfrac{1}{\beta}$-cocoercive w.r.t. $\|\cdot\|_S$ for $S\in\PD$ is also $\beta$-Lipschitz continuous w.r.t. $\|\cdot\|_S$.
  \label{prop:coco_to_Lip}
\end{prop}
\begin{pf}
  Follows immediately by definition of cocoercive operators and Cauchy-Schwarz inequality.
\end{pf}
The converse implication, however, does not hold in general.}
\begin{prop}
  An operator $T:\hilbert\to\hilbert$ that is $\tfrac{1}{\beta}$-cocoercive w.r.t. $\|\cdot\|_P$ for $P\in\PD$ is also $\beta\frac{\lambda_{\max}(S^{-1})}{\lambda_{\min}(P^{-1})}$-Lipschitz continuous w.r.t. $\|\cdot\|_S$ for all $S\in\PD$ and $\beta$-Lipschitz continuous w.r.t. $\|\cdot\|_P$.
  \label{prop:coco_to_Lip}
\end{prop}
\begin{pf}
  By $\tfrac{1}{\beta}$-cocoercivity of $T$ w.r.t. $\|\cdot\|_P$, we conclude:
  \begin{align*}
    \|Tx-Ty\|_{S^{-1}}^2&\leq\lambda_{\max}(S^{-1})\|Tx-Ty\|^2\leq\tfrac{\lambda_{\max}(S^{-1})}{\lambda_{\min}(P^{-1})}\|Tx-Ty\|_{P^{-1}}^2\\
                        &\leq\beta\tfrac{\lambda_{\max}(S^{-1})}{\lambda_{\min}(P^{-1})}\langle Tx-Ty,x-y\rangle\\
                        &\leq\beta\tfrac{\lambda_{\max}(S^{-1})}{\lambda_{\min}(P^{-1})}\|x-y\|_S\|Tx-Ty\|_{S^{-1}}.
  \end{align*}
The first result follows after rearrangement. The second follows trivially by Cauchy-Schwarz on the cocoercivity definition.
\end{pf}
The converse implication, that Lipschitz continuity implies cocoercivity, does not hold in general.
\begin{prop}
  Let $T:\hilbert\to\hilbert$ be $L$-Lipschitz continuous w.r.t. $\|\cdot\|$ and let $P,S\in\PD$ be arbitrary. Then
  \begin{align*}
    \|Tx-Ty\|_S\leq L\tfrac{\sqrt{\lambda_{\max}(S)}}{\sqrt{\lambda_{\min}(P)}}\|x-y\|_P
  \end{align*}
  for all $x,y\in\hilbert$.
  \label{prop:Lip}
\end{prop}
\begin{pf}
  For arbitrary $x,y\in\hilbert$, we conclude
  \begin{align*}
    \|Tx-Ty\|_S&\leq\sqrt{\lambda_{\max}(S)}\|Tx-Ty\|\leq L\sqrt{\lambda_{\max}(S)}\|x-y\|
               \leq L\tfrac{\sqrt{\lambda_{\max}(S)}}{\sqrt{\lambda_{\min}(P)}}\|x-y\|_P.
  \end{align*}
 \end{pf}

\begin{prop}
  Let $T:\hilbert\to\hilbert$ be $\sigma$-strongly monotone w.r.t. $\|\cdot\|_P$ and let $P,S\in\PD$ be arbitrary. Then
  \begin{align*}
    \|Tx-Ty\|_{S^{-1}}\geq\tfrac{\sigma\lambda_{\min}(P)}{\lambda_{\max}(S)}\|x-y\|_S
  \end{align*}
  for all $x,y\in\hilbert$.
  \label{prop:sm_lower}
\end{prop}
\begin{pf}
  By $\sigma$-strong monotonicity of $T$ w.r.t. $\|\cdot\|_P$, we conclude:
  \begin{align*}
    \sigma\|x-y\|_S^2&\leq\sigma\lambda_{\max}(S)\|x-y\|^2\leq\tfrac{\sigma\lambda_{\max}(S)}{\lambda_{\min}(P)}\|x-y\|_P^2\leq\tfrac{\lambda_{\max}(S)}{\lambda_{\min}(P)}\langle x-y,Tx-Ty\rangle\\
                     &\leq\tfrac{\lambda_{\max}(S)}{\lambda_{\min}(P)}\|x-y\|_{S}\|Tx-Ty\|_{S^{-1}}.
  \end{align*}
    Rearrangement this gives the result.
\end{pf}

\section{Nonlinear Forward-Backward Splitting}
\label{sec:NOFOB}

We propose nonlinear forward-backward splitting in Algorithm~\ref{alg:NOFOB} for solving monotone inclusions of the form
\begin{align}
  0\in Ax+Cx,
  \label{eq:prob}
\end{align}
that satisfy the following assumption.
\begin{ass}
Assume that:
\begin{enumerate}[(i)]
\item $A:\hilbert\to 2^\hilbert$ is maximally monotone.\label{ass:prob:A}
\item $C:\hilbert\to\hilbert$ is $\tfrac{1}{\beta}$-cocoercive w.r.t. $\|\cdot\|_P$ 
      with $\beta\in[0,4)$ for some $P\in\PD$.\label{ass:prob:C}
    
    \item The solution set $\zero(A+C):=\{x:0\in Ax+Cx\}$ is nonempty.\label{ass:prob:nonempty}
  \end{enumerate}
\label{ass:prob}
\end{ass}

The cocoercive term $C$ is assumed $\tfrac{1}{\beta}$-cocoercive w.r.t. $\|\cdot\|_P$, where $P$ is a self-adjoint positive definite operator that is part of the algorithm.
The inverse cocoercivity constant is constrained to satisfy $\beta\in[0,4)$. This is similar to the construction in AFBA \cite{Latafat2017Asymmetric} and is indeed no restriction. The choice of $P$ in the metric $\|\cdot\|_P$ is free and can always be chosen large enough to satisfy this.
Besides $P$, the algorithm also uses a possibly nonlinear iteration dependent kernel $M_k:\hilbert\to\hilbert$ for the forward-backward step and a linear operator $S\in\PD$ to define the projection metric. The algorithm is stated below. 
\ifx\one\showalgs\begin{algorithm}[H]
  \begin{algorithmic}[1]
    \STATE {\bf Input:} $x_0\in\hilbert$
		\FOR{$k=0,1,\ldots$}
                \STATE $\hat{x}_{k}:=(M_k+A)^{-1}(M_k-C)x_{k}$ 
                \STATE $H_k := \{z:\langle M_kx_k-M_k\hat{x}_k,z-\hat{x}_k\rangle\leq \tfrac{\beta}{4}\|x_k-\hat{x}_k\|_P^2\}$                
                \STATE $x_{k+1}:=  (1-\theta_k)x_k+\theta_k \proj_{H_k}^S(x_k)$\label{alg-step:rel_proj}
		\ENDFOR
	\end{algorithmic}
        \caption{Nonlinear Forward-Backward Splitting (NOFOB)}
        \label{alg:NOFOB}
\end{algorithm}\fi
The first step in the algorithm is a forward-backward step with a possibly nonlinear kernel $M_k$. If $C=0$, it is a pure backward step, that we also call a {\emph{nonlinear resolvent}} step. We will see that the nonlinear forward-backward step creates a (strictly, if $x_k\not\in\zero(A+C)$) separating hyperplane between the current point $x_k$ and the solution set $\zero(A+C)\neq\emptyset$. In the second step, the halfspace $H_k$ is constructed from the separating hyperplane. It contains the solution set but not $x_k$ and is called a separating halfspace. The halfspace construction requires an extra evaluation of $M_k$ at the point $\hat{x}_k$. The subsequent iterate $x_{k+1}$ is in the third step obtained by a relaxed projection from $x_k$ onto the created halfspace $H_k$. The relaxation is decided by $\theta_k\in(0,2)$ and $S$ defines the projection metric. 

We need the following assumptions on the kernels $M_k$ and the projection metric $\|\cdot\|_S$ for our convergence analysis. 
\begin{ass}
  Assume that all $M_k:\hilbert\to\hilbert$ are maximally monotone, single-valued, 1-strongly monotone w.r.t. $\|\cdot\|_P$ for some $P\in\PD$, and $L_M$-Lipschitz continuous (w.r.t. $\|\cdot\|$). Further assume that the projection metric $\|\cdot\|_S$ is defined by some $S\in\PD$.
\label{ass:alg}
\end{ass}

Note that the metric $\|\cdot\|_P$ is used both in the cocoercivity assumption of $C$ in Assumption~\ref{ass:prob} and in the strong monotonicity assumption of $M_k$. This restricts the choice of kernels $M_k$. In absence of $C$, there is no such restriction, and weaker assumptions suffice (see Remark~\ref{rem:M_not_strongly_mono}). 

That we allow for nonlinear and nonsymmetric kernels $M_k$ in the resolvents distinguishes our method from other methods based on forward-backward splitting. This is key for the construction of our four-operator splitting method that has, e.g., forward-backward-forward splitting~\cite{Tseng2000AModified}, forward-backward-half-forward splitting~\cite{BricenoArias2018Forward}, asymmetric forward-backward-adjoint splitting~\cite{Latafat2017Asymmetric,Latafat2018PrimalDual}, and synchronous projective splitting~\cite{Combettes2018Asynchronous} as special cases. 

\subsection{NOFOB with Explicit Projection}

\label{sec:NOFOB_explicit}

We can derive an explicit formula for the projection step in Algorithm~\ref{alg:NOFOB}. The result of the projection is found by solving
\begin{align*}
    \begin{tabular}{ll}
        $\displaystyle\minimize_z$ & $\|z-x_k\|_S^2$\\
        subject to &$\langle M_kx_k-M_k\hat{x}_k,z-\hat{x}_k\rangle\leq \tfrac{\beta}{4}\|x_k-\hat{x}_k\|_P^2$.
    \end{tabular}
\end{align*}
Assuming that $x_k\not\in H_k$, the solution $z$ can be found by projection onto the separating hyperplane, i.e.,
\begin{align}
  z&=x_k-\frac{\langle M_kx_k-M_k\hat{x}_k,x_k-\hat{x}_k\rangle-\tfrac{\beta}{4}\|x_k-\hat{x}_k\|_P^2}{\|M_kx_k-M_k\hat{x}_k\|_{ S^{-1}}^2}S^{-1}(M_kx_k-M_k\hat{x}_k).
     \label{eq:proj_formula}
\end{align}
Therefore, by letting
\begin{align}
  \mu_k:=\frac{\langle M_kx_k-M_k\hat{x}_k,x_k-\hat{x}_k\rangle-\tfrac{\beta}{4}\|x_k-\hat{x}_k\|_P^2}{\|M_kx_k-M_k\hat{x}_k\|_{ S^{-1}}^2},
  \label{eq:muk}
\end{align}
the projection step can be written explicitly as in the following algorithm.
\ifx\one\showalgs\begin{algorithm}[H]
  \begin{algorithmic}[1]
    \STATE {\bf Input:} $x_0\in\hilbert$
		\FOR{$k=0,1,\ldots$}
                \STATE $\hat{x}_{k}:=(M_k+A)^{-1}(M_k-C)x_{k}$ 
                \STATE $\mu_k:=\frac{\langle M_kx_k-M_k\hat{x}_k,x_k-\hat{x}_k\rangle-\tfrac{\beta}{4}\|x_k-\hat{x}_k\|_P^2}{\|M_kx_k-M_k\hat{x}_k\|_{ S^{-1}}^2}$
                \STATE $x_{k+1}:=x_k-\theta_k\mu_kS^{-1}(M_kx_k-M_k\hat{x}_k)$
		\ENDFOR
	\end{algorithmic}
        \caption{NOFOB with Explicit Projection}
        \label{alg:NOFOB_explicit}
      \end{algorithm}\fi

      
\subsection{NOFOB with Conservative Step}

\label{sec:NOFOB_CS}

Algorithm~\ref{alg:NOFOB_explicit} will converge if $\mu_k$ defined in \eqref{eq:muk} is replaced by $\hat{\mu}_k\geq\epsilon_{\mu}>0$ that satisfies
\begin{align}
  \hat{\mu}_k\leq\frac{\langle M_kx-M_ky,x-y\rangle-\tfrac{\beta}{4}\|x-y\|_P^2}{\|M_kx-M_ky\|_{ S^{-1}}^2}
  \label{eq:mu_bound}
\end{align}
for all $x,y\in\hilbert$ such that $x\neq y$. Unless $k\in\naturals$ is such that $x_k=\hat{x}_k$, we immediately see that $\hat{\mu}_k\leq\mu_k$ with $\mu_k$ in \eqref{eq:muk}. For $k\in\naturals$ such that $x_k=\hat{x}_k$, the choice of real-valued $\hat{\mu}_k$ does not matter since it is multiplied by 0 in Algorithm~\ref{alg:NOFOB_explicit} so we can still use a $\hat{\mu}_k$ that satisfies \eqref{eq:mu_bound}. As we will see in Proposition~\ref{prop:muk_bounds}, a uniform lower bound $\epsilon_{\mu}>0$ can be chosen such that $\hat{\mu}_k\geq \epsilon_{\mu}$. Therefore, we can replace $\mu_k$ in Algorithm~\ref{alg:NOFOB_explicit} by $\hat{\mu}_k$, which leads to the following algorithm.

\ifx\one\showalgs\begin{algorithm}[H]
  \begin{algorithmic}[1]
    \STATE {\bf Input:} $x_0\in\hilbert$, $\hat{\mu}_k>\epsilon_{\mu}>0$ such that \eqref{eq:mu_bound} holds for all $k\in\naturals$ and all $x\neq y\in\hilbert$
    \FOR{$k=0,1,\ldots$}
                \STATE $\hat{x}_{k}:=(M_k+A)^{-1}(M_k-C)x_{k}$ 
                \STATE $x_{k+1}:=x_k-\theta_k\hat{\mu}_k S^{-1}(M_kx_k-M_k\hat{x}_k)$
		\ENDFOR
	\end{algorithmic}
        \caption{NOFOB with Conservative Step}
        \label{alg:NOFOB_constant_step}
      \end{algorithm}\fi

      The only effective change in the algorithm compared to Algorithm~\ref{alg:NOFOB_explicit} is that a smaller relaxation factor is used for the projection. 
      If $x_{k}\neq\hat{x}_k$, we know that $\mu_k>0$ (Proposition~\ref{prop:muk_bounds}) and the last update step can be written as
\begin{align*}
  x_{k+1}=x_k-\theta_k\hat{\mu}_kS^{-1}(M_kx_k-M_k\hat{x}_k)=x_k-\tfrac{\theta_k\hat{\mu}_k}{\mu_k}\mu_k S^{-1}(M_kx_k-M_k\hat{x}_k).
\end{align*}
This reveals that Algorithm~\ref{alg:NOFOB_constant_step} is equivalent to Algorithm~\ref{alg:NOFOB_explicit} with the smaller relaxation parameter $\tfrac{\theta_k\hat{\mu}_k}{\mu_k}$. The relaxation parameter is directly proportional to the step-length taken by the algorithm, so Algorithm~\ref{alg:NOFOB_constant_step} is a short-step version of Algorithm~\ref{alg:NOFOB_explicit}, where the reduction in step-length for every $k\in\naturals$ is exactly the ratio $\tfrac{\hat{\mu}_k}{\mu_k}\leq 1$. The $\mu_k$ in \eqref{eq:muk} is a {\emph{local}} property parameter of $M_k$ that needs to hold only for $x_k$ and $\hat{x}_k$. The parameter $\hat{\mu}_k$ in \eqref{eq:mu_bound} is, for every $k$, a {\emph{global}} parameter for the same property of $M_k$ that must hold for all $x\neq y\in\hilbert$. The ratio between the global and local property parameters decide the reduction in step-length, which can be significant. As we will see, standard formulations of FBF \cite{Tseng2000AModified} and FBHF \cite{BricenoArias2018Forward} are  special cases of Algorithm~\ref{alg:NOFOB_constant_step}, so they are conservative special cases of Algorithm~\ref{alg:NOFOB_explicit}.

\section{Nonlinear Forward-Backward Map}

\label{sec:FBmap}

In this section, we analyze properties of the nonlinear forward-backward map 
\begin{align}
  \FB:=(M+A)^{-1}\circ (M-C)
  \label{eq:FB}
\end{align}
that is used in the first step in NOFOB and generalizes the so-called $F$-resolvent in \cite{Bauschke2010General} to include a cocoercive term $C$. We analyze the map for kernel $M$ without iteration index for cleaner notation. We show that $\FB$ has full domain, is Lipschitz continuous, and is single-valued. We also state the straightforward result that the fixed-point set of $\FB$ coincides with the zero set of $A+C$. Finally, we show that application of $\FB$ at $x$ creates a strictly separating hyperplane between $x$ and $\zero(A+C)$ for all $x\not\in\zero(A+C)$.
%
\begin{prop}
  Suppose that Assumptions~\ref{ass:prob}~and~\ref{ass:alg} hold. Then $\FB$:
  \begin{enumerate}[(i)]
  \item has full domain,\label{item:full_dom}
  \item is Lipschitz continuous,\label{item:Lipschitz}
  \item is single-valued.\label{item:single_valued}
  \end{enumerate}
\label{prop:FB_prop}
\end{prop}
\begin{pf}
  The point $\hat{x}:=\FB x$ is obtained by applying $(M+A)^{-1}$ to $(M-C)x$. Since $M$ and $C$ by assumption have full domains, are Lipschitz continuous, and single-valued, it suffices to show that $(M+A)^{-1}$ satisfies {\it(\ref{item:full_dom})-(\ref{item:single_valued})}.

  We use that 1-strong monotonicity of $M$ w.r.t. $\|\cdot\|_P$ implies strong monotonicity w.r.t. $\|\cdot\|$ of $M$ and of $A+M$. Further, since $M$ has full domain and $M$ and $A$ are maximally monotone, so is $M+A$ \cite[Corollary~25.5]{Bauschke2017Convex}. Therefore, $M+A$ is maximally strongly monotone.
  
  {\it(\ref{item:full_dom})} Applying \cite[Corollary~25.28]{Bauschke2017Convex} implies that $\dom(M+A)^{-1}=\ran(M+A)=\hilbert$.

{\it(\ref{item:Lipschitz})}: Applying \cite[Example~22.7]{Bauschke2017Convex} gives that $(M+A)^{-1}$ is cocoercive, hence Lipschitz continuous (Proposition~\ref{prop:coco_to_Lip}). 

  {\it(\ref{item:single_valued})} Follows from {\it(\ref{item:full_dom})} and {\it(\ref{item:Lipschitz})}.


\end{pf}
The Lipschitz constant of $\FB$ is typically larger than one. Hence $\FB$ is not nonexpansive in general. This is in contrast to the standard forward-backward map that is nonexpansive (averaged) for appropriate step-sizes.

A consequence of Proposition~\ref{prop:FB_prop} is that NOFOB will generate an infinite sequence $(x_k)_{k\in\naturals}$. The output of the nonlinear forward-backward step $\hat{x}_k$ exists ($\FB$ has full domain) and is unique ($\FB$ is single-valued). The subsequent projection onto the separating halfspace also exists and is unique. Therefore, the algorithm will not come to a halt, but instead generate an infinite sequence of points $(x_k)_{k\in\naturals}$.

Next, we show that the fixed-point set of $\FB$ coincides with the zero set of $A+C$.
\begin{prop}
  Suppose that Assumptions~\ref{ass:prob}~and~\ref{ass:alg} hold. Then
  \begin{align*}
    \fix \FB=\zero (A+C).
  \end{align*}
  \label{prop:fixedFB-zeroAC}
\end{prop}
\begin{pf}
    We know from Proposition~\ref{prop:FB_prop} that $\FB=(M+A)^{-1}\circ(M-C)$ is single-valued and has full domain. Hence, $x\in\fix \FB$ is equivalent to
    \begin{align*}
        x=(M+A)^{-1}(M-C)x\Leftrightarrow Mx+Ax\ni Mx-Cx \Leftrightarrow 0\in Ax+Cx,
    \end{align*}
    since $M$ and $C$ are single-valued. This concludes the proof.
\end{pf}

Next, we show that the forward-backward map creates a separating hyperplane between the point of application and the solution set. The hyperplane is for every $x\in\hilbert$ defined as the 0\ordth level of the affine function
\begin{align}
  \psi_{x}(z):=\langle Mx-M\FB x,z-\FB x\rangle-\tfrac{\beta}{4}\|x-\FB x\|_P^2,
  \label{eq:psi_def_FB}
\end{align}
where $\beta\in[0,4)$ is the cocoercivity constant of $C$ in Assumption~\ref{ass:prob}.

\begin{thm}
Suppose that Assumptions~\ref{ass:prob}~and~\ref{ass:alg} hold and let $\psi_x$ be as in \eqref{eq:psi_def_FB}. Then 
\begin{enumerate}[(i)]
\item $\psi_x(x)\geq(1-\tfrac{\beta}{4})\|x-\FB x\|_{P}^2\geq 0$, i.e., $\psi_x(x)\geq 0$, for all $x\in\hilbert$,
\item $\psi_x(x)=0$ if and only if $x\in\zero(A+C)=\fix \FB$,
\item $\psi_x(z)\leq 0$ for all $z\in\zero(A+C)=\fix \FB$ and $x\in\hilbert$.\label{item:psi_negative_solutions}
\end{enumerate}
\label{thm:sep-hyper}
\end{thm}
\begin{pf}
  {\it(i)}
  It holds by 1-strong monotonicity of $M$ w.r.t. $\|\cdot\|_P$ (Assumption~\ref{ass:alg}) that
  \begin{align*}
    \psi_x(x)&=\langle Mx-M\FB x,x-\FB x\rangle-\tfrac{\beta}{4}\|x-\FB x\|_P^2
               \geq (1-\tfrac{\beta}{4})\|x-\FB x\|_P^2,
  \end{align*}
  Now, since $\beta\in[0,4)$ this quantity is nonnegative, which proves {\it(i)}.

  {\it(ii)} That $\psi_x(x)>0$ for all $x\not\in\fix\FB$ follows from {\it(i)} and the definition of a fixed-point. Insertion of any fixed-point $z=x=\FB x$ in \eqref{eq:psi_def_FB} gives $\psi_x(x)=0$. An appeal to Proposition~\ref{prop:fixedFB-zeroAC} proves {\it(ii)}.

  {\it(iii)}  Since $\FB$ is single-valued and has full domain (Proposition~\ref{prop:FB_prop}), it holds by the definition of $\FB$ and since $M$ is single-valued, that $M\FB x+A\FB x\ni Mx-Cx$, or equivalently $A\FB x\ni Mx-Cx-M\FB x$, for all $x\in\hilbert$. Since $z\in\fix \FB=\zero(A+C)\neq\emptyset$ we have $-Cz\in Az$. By monotonicity of $A$ evaluated at $z$ and $\FB x$, we therefore get
  \begin{align}
    \nonumber 0&\geq-\langle -Cz-(Mx-M\FB x+Cx),z-\FB x\rangle\\
     &=\langle Mx-M\FB x,z-\FB x\rangle+\langle Cz-Cx,z-\FB x\rangle.
       \label{eq:AmonoIneq}
  \end{align}
  Now, if $\beta=0$, $C$ is constant, the second term is 0, and the result holds. For $\beta\in(0,4)$, we have 
  \begin{align*}
    \langle Cz-Cx,z-\FB x\rangle
    &=\langle Cz-Cx,x-\FB x\rangle+\langle Cz-Cx,z-x\rangle\\
    &\geq -\tfrac{\epsilon}{2}\|x-\FB x\|_{P}^2-\tfrac{1}{2\epsilon}\|Cx-Cz\|_{P^{-1}}^2+\langle Cx-Cz,x-z\rangle\\
    &\geq -\tfrac{\epsilon}{2}\|x-\FB x\|_{P}^2-(\tfrac{1}{2\epsilon}-\tfrac{1}{\beta})\|Cx-Cz\|_{P^{-1}}^2\\
    &=-\tfrac{\beta}{4}\|x-\FB x\|_{P}^2,
  \end{align*}
where we have used negated Young's inequality, $\tfrac{1}{\beta}$-cocoercivity of $C$ w.r.t. $\|\cdot\|_P$, and $\epsilon=\tfrac{\beta}{2}>0$. Inserting into \eqref{eq:AmonoIneq} gives the result.
\end{pf}

The halfspace $H_k$ that is projected onto in Algorithm~\ref{alg:NOFOB} can be written as
\begin{align}
  H_k = \{z:\psi_{x_k}(z)\leq 0\},
  \label{eq:Hk_def}
\end{align}
with iteration dependent kernel $M_k$ in the definitions of $\FB$ and $\psi_x$ in \eqref{eq:psi_def_FB}. Theorem~\ref{thm:sep-hyper} shows that $H_k$ is a separating halfspace (i.e., constructed from a separating hyperplane) between the current iterate $x_k$ and the solution set $\zero(A+C)$. The separation is strict unless $x_k$ has already solved the inclusion problem, i.e., unless $x_k\in\zero(A+C)$. Separation is the key property that allows us to show convergence of Algorithm~\ref{alg:NOFOB} in Section~\ref{sec:convergence} by means of \Fejer~monotonicity.

\begin{rem}
  \label{rem:M_not_strongly_mono}
  The requirement that $M$ is strongly monotone w.r.t. $\|\cdot\|_P$ can in some cases be relaxed without loosing strict separation. It is evident from the proof of Theorem~\ref{thm:sep-hyper} that $\psi_x(z)\leq 0$ for all $z\in\zero(A+C)$ under very mild assumptions on $M$. We do not even need monotonicity. Strong monotonicity is needed to ensure that $\psi_x(x)\geq(1-\tfrac{\beta}{4})\|x-\FB x\|_P^2>0$ when $\beta\in[0,4)$ to have separation. If no cocoercive term is present, i.e., $C=0$, then $\beta=0$ and strict monotonicity of $M$ would give
  \begin{align}
    \psi_x(x)=\langle Mx-M\FB x,x-\FB x\rangle>0,
    \label{eq:sep_strictly_mono}
  \end{align}
  which shows strict separation. To guarantee full domain and single-valuedness of the nonlinear resolvent map $(M+A)^{-1}M$, we can, e.g., further assume that $M$ is $3^*$-monotone and surjective \cite[Proposition~4.2]{Bauschke2010General}. These assumptions combined are weaker than strong monotonicity w.r.t. $\|\cdot\|_P$. The stronger assumption is used here to have separation in presence of cocoercive term $C$.
\end{rem}
\begin{rem}
  The setting with $C=0$ is also related to Bregman forward-backward steps that are used, e.g., in the algorithms in \cite{Bauschke2017ADescent,Lu2018Relatively}. These works consider 
  \begin{align*}
    \minimize_{x\in C} f(x)+g(x),
  \end{align*}
  where $f:\hilbert\to\reals\cup\{\infty\}$ and $g:\mathcal{D}\to\hilbert$ are proper closed and convex and $C\subset \hilbert$ is a nonempty closed convex set. Let us assume that the solution set is nonempty, that $C=\overline{\dom} h$ for a strictly convex function $h$, and that $g$ is continuously differentiable on $\interior\,C$ with $\mathcal{D}=\dom g\supseteq \dom h$ and $\beta$-smooth w.r.t. $h$, which is defined as that $\beta h-g$ is convex on $\interior\,C$ (see \cite{Bauschke2017ADescent,Lu2018Relatively}) or equivalently that $\beta\nabla h-\nabla g$ is monotone on $\interior\,C$. Letting $F=\gamma^{-1}h-g$, $M=\nabla F=\gamma^{-1}\nabla h-\nabla g$, and $A=\partial f+\nabla g$ gives that the resolvent step
  \begin{align*}
    \hat{x}&=(M+A)^{-1}Mx=(\nabla F+\partial f+\nabla g)^{-1}\nabla Fx
           =(\gamma^{-1}\nabla h+\partial f)^{-1}(\gamma^{-1}\nabla h-\nabla g)x\\
           &=(\nabla h+\gamma\partial f)^{-1}(\nabla h-\gamma\nabla g)x,
  \end{align*}
  which is a Bregman forward-backward step with step-size $\gamma$. For $\gamma\in(0,\tfrac{1}{\beta})$, we have that $M=\nabla F=\gamma^{-1}\nabla h-\nabla g=(\gamma^{-1}-\beta)\nabla h+(\beta\nabla h-\nabla g)$, which is strictly monotone since $(\gamma^{-1}-\beta)>0$, $\nabla h$ is strictly monotone, and $\beta\nabla h-\nabla g$ is monotone by assumption. We conclude from \eqref{eq:sep_strictly_mono} and Theorem~\ref{thm:sep-hyper}({\it\ref{item:psi_negative_solutions}}), which can be shown to hold for all $x\in\dom\FB$, that we have strict separation. Although we will not cover algorithms for this setting, it is interesting to note that the Bregman forward-backward step considered in \cite{Bauschke2017ADescent,Lu2018Relatively} can be interpreted as a pure Bregman backward step with $F=\gamma^{-1}h-g$, which is a special case of the nonlinear resolvent.
\end{rem}
\begin{rem}
An interpretation of the nonlinear forward-backward map arises by letting $E:=M-P$. Then $E$ is single-valued and maximal monotone since $M$ is assumed maximally monotone and strongly monotone w.r.t. $\|\cdot\|_P$. We can write the nonlinear forward-backward map applied to $A$ and $C$ as a standard forward-backward map with linear positive definite kernel $P$ applied to the maximal monotone $(A+E)$ and single-valued, but typically not cocoercive, $C-E$:
\begin{align*}
  \hat{x}=(M+A)^{-1}(M-C)x=(P+(A+E))^{-1}(P-(C-E))x.
\end{align*}
Although interesting, this interpretation has some limitations and we cannot use it to simplify proofs, which is why we will not utilize it further. First, it does not in a natural way allow for the kernel $M$ to be iteration dependent in algorithms. This interpretation would give a {\emph{problem}} that is dependent on the algorithm iterate. Second, since $C-E$ is not cocoercive, we cannot make use of a standard forward-backward analysis that is based on averaged operators, so proofs will be the same. Third, when $C=0$, we no longer need strong monotonicity of $M$, see Remark~\ref{rem:M_not_strongly_mono}. Without strong monotonicity, the interpretation is no longer valid.
\end{rem}


\section{Convergence}

\label{sec:convergence}


We first present a general result on weak convergence of sequences to the solution set $\zero(A+C)$. We then show that Algorithm~\ref{alg:NOFOB} satisfies the required assumptions. Before we state the result, we introduce the iteration dependent nonlinear forward-backward map
\begin{align}
  \FB^{k}:=(M_k+A)^{-1}(M_k-C),
  \label{eq:FBk}
\end{align}
that is used in the first step in NOFOB.
\begin{prop}
  Suppose that Assumptions~\ref{ass:prob}~and~\ref{ass:alg} hold. Assume that $(x_k)_{k\in\naturals}$ and $(\hat{x}_k)_{k\in\naturals}$ are sequences in $\hilbert$ satisfying $\hat{x}_k=\FB^kx_k$, and that
  \begin{enumerate}[(i)]
  \item $(\|x_k-z\|_S)_{k\in\naturals}$ converges for every $z\in\zero(A+C)$,\label{item:dist_conv}
  \item $(\|x_k-\hat{x}_k\|_{S})_{k\in\naturals}$ converges to 0. \label{item:resid_conv}
  \end{enumerate}
   Then $x_k\weakto \bar{x}\in\zero(A+C)$.
   \label{prop:seq_conv}
 \end{prop}
 \begin{pf}
   By 
   \cite[Theorem~4.11 and Example~4.6]{Bauschke2003Bregman}, it is enough to show that every weak sequential cluster point belongs to $\zero(A+C)$. (These results in \cite{Bauschke2003Bregman} are actually much more general than what we need. Another way to arrive at this conclusion is to generalize \cite[Lemma~2.47]{Bauschke2017Convex} to handle also scaled norms $\|\cdot\|_S$ as in {\it (i)}. This generalization is straightforward.)

   At least one cluster point exists since {\it (i)} and Assumption~\ref{ass:prob} imply that $(x_k)_{k\in\naturals}$ is bounded. Let $(x_{n_k})_{k\in\naturals}$ be a weakly convergent subsequence $x_{n_k}\weakto x$, where $x$ is the cluster point. Since $\hat{x}_{n_k}=\FB^{n_k} x_{n_k}=(M_{n_k}+A)^{-1}(M_{n_k}-C)x_{n_k}$, it holds that $(M_{n_k}+A)\hat{x}_{n_k}\ni (M_{n_k}-C)x_{n_k}$ by single-valuedness of $M_{n_k}$ and $C$. This implies that 
   \begin{align}
     u_{n_k}:=M_{n_k}x_{n_k}-M_{n_k}\hat{x}_{n_k}-(Cx_{n_k}-C\hat{x}_{n_k})\in (A+C)\hat{x}_{n_k}.
     \label{eq:inGraph}
   \end{align}
   Since $x_{n_k}\weakto x$ and, by item~{\it (ii)}, $x_{n_k}-\hat{x}_{n_k}\to 0$, we conclude $\hat{x}_{n_k}\weakto x$. Further, by $L_M$-Lipschitz continuity w.r.t. $\|\cdot\|$ of all $M_k$ (Assumption~\ref{ass:alg}) and $\tfrac{\beta}{\lambda_{\min}(P^{-1})}$-Lipschitz continuity w.r.t. $\|\cdot\|$ of $C$ (Assumption~\ref{ass:prob} and Proposition~\ref{prop:coco_to_Lip}) and since $\|x_{n_k}-\hat{x}_{n_k}\|\to 0$, we conclude that
   \begin{align*}
     \|u_{n_k}\|&=\|M_{n_k}x_{n_k}-M_{n_k}\hat{x}_{n_k}-(Cx_{n_k}-C\hat{x}_{n_k})\|\\
                &\leq \|M_{n_k}x_{n_k}-M_{n_k}\hat{x}_{n_k}\|+\|Cx_{n_k}-C\hat{x}_{n_k}\|\\
                &\leq \left(L_M+\tfrac{\beta}{\lambda_{\min}(P^{-1})}\right)\|x_{n_k}-\hat{x}_{n_k}\|\to 0,
   \end{align*}
   i.e., $u_{n_k}\to 0$. Now, since $A+C$ is maximally monotone (since $C$ has full domain, see \cite[Corollary~25.5]{Bauschke2017Convex}) and $(\hat{x}_{n_k},u_{n_k})\in\gph(A+C)$, we conclude by weak-strong closedness of graphs of maximally monotone operators \cite[Proposition~20.38]{Bauschke2017Convex} that $(x,0)\in\gph(A+C)$, i.e., $x\in\zero(A+C)$. That is, every weak sequential cluster point belongs to $\zero(A+C)$ and the proof is complete.
 \end{pf}

 Before we prove convergence of Algorithm~\ref{alg:NOFOB}, we need the following bounds on $\mu_k$ defined in \eqref{eq:muk}. The proof is deferred to Appendix~\ref{app:muk_bounds_pf}.
 
 \begin{prop}
   Suppose that Assumptions~\ref{ass:prob}~and~\ref{ass:alg} hold. Then for all $x,y\in\hilbert$ such that $x\neq y$,
   \begin{align*}
     \mu_k^{x,y}:=\frac{\langle M_kx-M_ky,x-y\rangle-\tfrac{\beta}{4}\|x-y\|_P^2}{\|M_kx-M_ky\|_{S^{-1}}^2}
   \end{align*}
   satisfies 
   \begin{align*}
     \mu_k^{x,y}\in\left[(1-\tfrac{\beta}{4})\tfrac{\lambda_{\min}(P)}{L_M^2\lambda_{\max}(S^{-1})},\tfrac{\lambda_{\max}(S)}{\lambda_{\min}(P)}\right].
   \end{align*}
In particular, the dual variable in the $H_k$ projection, $\mu_k=\mu_{k}^{x_k,\hat{x}_k}$ in \eqref{eq:muk}, satisfies the same bounds 
   for all $x_k\not\in\zero(A+C)$.
   \label{prop:muk_bounds}
 \end{prop}


 \begin{thm}
   Suppose that Assumptions~\ref{ass:prob}~and~\ref{ass:alg} hold and that the relaxation parameter $\theta_k\in(0,2)$ satisfies $\liminf_{k\to\infty}\theta_k(2-\theta_k)>0$. Then Algorithm~\ref{alg:NOFOB} constructs sequences $(x_k)_{k\in\naturals}$ and $(\hat{x}_k)_{k\in\naturals}$ such that
   \begin{enumerate}[(i)]
   \item $(\|x_k-z\|_S)_{k\in\naturals}$ converges for every $z\in\zero(A+C)$,
   \item $(\|x_k-\hat{x}_k\|_S)_{k\in\naturals}$ converges to 0,
   \item $(x_k)_{k\in\naturals}$ converges weakly to a point in $\zero(A+C)$, i.e., $x_k\weakto \bar{x}\in\zero(A+C)$ as $k\to\infty$.
     \end{enumerate}
   \label{thm:conv}
\end{thm}
\begin{pf}
Proposition~\ref{prop:FB_prop} guarantees that infinite sequences $(x_k)_{k\in\naturals}$ and $(\hat{x}_k)_{k\in\naturals}$ are constructed by Algorithm~\ref{alg:NOFOB}. First assume that there exits $k\in\naturals$ such that $x_k=\bar{x}\in\zero(A+C)$. Then $\hat{x}_k=x_k\in\zero(A+C)\subseteq H_k$ (Theorem~\ref{thm:sep-hyper}) and therefore $x_{k+1}=x_k=\bar{x}\in\zero(A+C)$. Induction gives that the sequence will stay at $\bar{x}$, and that {\it(i)}-{\it(iii)} hold. It is left to prove weak convergence when all iterates $x_k\not\in\zero(A+C)$.

{\it(i)}: Let us define the relaxed projection operator $\Pi_{H_k,\theta_k}^S:=(1-\theta_k)\id+\theta_k\Pi_{H_k}^S$ that satisfies $x_{k+1}=\Pi_{H_k,\theta_k}^Sx_k$ and $z=\Pi_{H_k,\theta_k}^Sz$ for all $z\in H_k$ and note that $\Pi_{H_k,\theta_k}^Sx_k-x_k:=\theta_k(\Pi_{H_k}^Sx_k-x_k)$. With an appeal to \cite[Proposition~4.16, Corollary~4.41, Proposition~4.35]{Bauschke2017Convex} we conclude that the relaxed projection operator is $\tfrac{\theta_k}{2}$-averaged w.r.t. $\|\cdot\|_S$ and satisfies
  \begin{align}
    \|x_{k+1}-z\|_S^2\leq\|x_k-z\|_S^2-\theta_k(2-\theta_k)\|x_k-\Pi_{H_k}^Sx_{k}\|_S^2
    \label{eq:seq-Fejer}
  \end{align}
  for all $z\in H_k$. In particular, it holds for all $z\in\zero(A+C)$ since $\zero(A+C)\subseteq H_k$ (Theorem~\ref{thm:sep-hyper}). Hence, $(\|x_k-z\|_S^2)_{k\in\naturals}$ is nonincreasing. Since it is also lower bounded, it converges.
  
{\it(ii)}:
First note that \eqref{eq:seq-Fejer} implies summability of $(\theta_k(2-\theta_k)\|x_k-\Pi_{H_k}^Sx_{k}\|_S^2)_{k\in\naturals}$. Since the relaxation parameter $\theta_k\in(0,2)$ satisfies $\liminf_{k\to\infty}\theta_k(2-\theta_k)>0$, we conclude that $\|x_k-\Pi_{H_k}^Sx_k\|_S\to 0$. Now, since $x_k\not\in\zero(A+C)$, we have $x_k\not\in H_k$ (Theorem~\ref{thm:sep-hyper}), and the projection formula is given by \eqref{eq:proj_formula}. Therefore, using Proposition~\ref{prop:muk_bounds} and Proposition~\ref{prop:sm_lower}, we conclude
\begin{align}
  \nonumber 0\leftarrow\|x_k-\Pi_{H_k}^Sx_k\|_S&=\mu_k\|S^{-1}(M_kx_k-M_k\hat{x}_k)\|_S\\
  \nonumber                   &\geq(1-\tfrac{\beta}{4})\tfrac{\lambda_{\min}(P)}{L_M^2\lambda_{\max}(S^{-1})}\|M_kx_k-M_k\hat{x}_k\|_{S^{-1}}\\
                   &\geq(1-\tfrac{\beta}{4})\tfrac{\lambda_{\min}(P)}{L_M^2\lambda_{\max}(S^{-1})}\tfrac{\lambda_{\min}(P)}{\lambda_{\max}(S)}\|x_k-\hat{x}_k\|_{S}
\label{eq:proj_diff_to_residual}   
\end{align}
where the factor in front of $\|x_k-\hat{x}_k\|_S$ is positive since $\beta\in[0,4)$ by Assumption~\ref{ass:prob} and since $S,P\in\PD$. Hence $\|x_k-\hat{x}_k\|_S\to 0$.

{\it(iii)}: Follows from {\it(i)} and {\it(ii)} and Proposition~\ref{prop:seq_conv}.

\end{pf}

The following corollaries show convergence of Algorithm~\ref{alg:NOFOB_explicit} and its conservative variation in Algorithm~\ref{alg:NOFOB_constant_step}.

\begin{cor}
  Suppose that the assumptions of Theorem~\ref{thm:conv} hold. Then Algorithm~\ref{alg:NOFOB_explicit} is equivalent to Algorithm~\ref{alg:NOFOB} and Algorithm~\ref{alg:NOFOB_explicit} generates a sequence $(x_k)_{k\in\naturals}$ that converges weakly to a point in $\zero(A+C)$, i.e., $x_k\weakto \bar{x}\in\zero(A+C)$ as $k\to\infty$.
  \label{cor:conv_exp}
\end{cor}
\begin{pf}
  $x_k\not\in\zero(A+C)$: Theorem~\ref{thm:sep-hyper} implies that $x_k\not\in H_k$, and the projection formula \eqref{eq:proj_formula} used in Algorithm~\ref{alg:NOFOB_explicit} is exactly the projection used in Algorithm~\ref{alg:NOFOB}. 

  $x_k=\bar{x}\in\zero(A+C)$: Theorem~\ref{thm:sep-hyper} implies that $x_k\in H_k$ and $x_{k+1}=x_k$ in Algorithm~\ref{alg:NOFOB}. For Algorithm~\ref{alg:NOFOB_explicit}, we note that $x_k=\hat{x}_k$ (Proposition~\ref{prop:fixedFB-zeroAC}) and therefore $M_kx_k-M_k\hat{x}_k=0$, which implies that $x_{k+1}=x_k$, independent on $\hat{\mu}_k$. We conclude that the algorithms are equivalent and Theorem~\ref{thm:conv} shows convergence.
\end{pf}

 \begin{cor}
   Suppose that Assumptions~\ref{ass:prob}~and~\ref{ass:alg} hold, that the relaxation parameter $\theta_k\in(0,2)$ satisfies $\liminf_{k\to\infty}\theta_k(2-\theta_k)>0$, and that $\hat{\mu}_k$ in Algorithm~\ref{alg:NOFOB_constant_step} satisfies for all $k\in\naturals$ that $\hat{\mu}_k\geq\epsilon_{\mu}$ for some $\epsilon_{\mu}>0$ and that \eqref{eq:mu_bound} holds for all $x,y\in\hilbert$ such that $x\neq y$. Then Algorithm~\ref{alg:NOFOB_constant_step} generates a sequence $(x_k)_{k\in\naturals}$ that converges weakly to a point in $\zero(A+C)$, i.e., $x_k\weakto \bar{x}\in\zero(A+C)$ as $k\to\infty$.
   \label{cor:conv_conserv}
 \end{cor}
 \begin{pf}
   From Proposition~\ref{prop:muk_bounds} we conclude that there exists $\epsilon_{\mu}>0$ and $\hat{\mu}_k$ such that for all $k\in\naturals$, $\hat{\mu}_k\geq\epsilon_{\mu}$ and $\hat{\mu}_k$ satisfies \eqref{eq:mu_bound} for all $x,y\in\hilbert$ such that $x\neq y$. 
   
   Let $x_k\not\in\zero(A+C)$, then $\mu_k>0$ (Proposition~\ref{prop:muk_bounds}) and $\theta_{k}\hat{\mu}_k=\tfrac{\theta_{k}\hat{\mu}_k}{\mu_k}\mu_k$. Hence Algorithm~\ref{alg:NOFOB_constant_step} is Algorithm~\ref{alg:NOFOB_explicit}, and therefore Algorithm~\ref{alg:NOFOB} (see Corollary~\ref{cor:conv_exp}), with relaxation parameter $\tfrac{\theta_{k}\hat{\mu}_k}{\mu_k}$. Assuming that $x_k\not\in\zero(A+C)$ for all $k\in\naturals$, we show convergence via Theorem~\ref{thm:conv} by showing that $\liminf_{k\to\infty}\tfrac{\theta_{k}\hat{\mu}_k}{\mu_k}(2-\tfrac{\theta_{k}\hat{\mu}_k}{\mu_k})>0$. Let $\overline{\mu}>0$ be the finite upper bound for $\mu_k$ in Proposition~\ref{prop:muk_bounds}. Then 
   \begin{align*}
     \liminf_{k\to\infty}\tfrac{\theta_{k}\hat{\mu}_k}{\mu_k}(2-\tfrac{\theta_{k}\hat{\mu}_k}{\mu_k})\geq\tfrac{\epsilon_{\mu}}{\overline{\mu}}\liminf_{k\to\infty}\theta_{k}(2-\tfrac{\theta_{k}\hat{\mu}_k}{\mu_k})\geq\tfrac{\epsilon_{\mu}}{\overline{\mu}}\liminf_{k\to\infty}\theta_{k}(2-\theta_{k})>0.
   \end{align*}

   Let $x_k=\bar{x}\in\zero(A+C)$. Then $\hat{x}_k=x_k$ (Proposition~\ref{prop:fixedFB-zeroAC}) and therefore $x_{k+1}=x_k=\bar{x}$ and $x_{k+l}=\bar{x}$ for all $l\in\naturals$ and the convergence result trivially holds.
 \end{pf}

\subsection{Linear convergence}

In this section, we show local linear convergence of Algorithm~\ref{alg:NOFOB} under metric subregularity. We define metric subregularity as follows.
\begin{defin}[Metric subregularity]
  A mapping $F:\hilbert\to 2^\hilbert$ is metrically subregular at $\bar{x}$ for $\bar{y}$ if $\bar{y}\in F\bar{x}$ and there exists $\kappa\in[0,\infty)$ along with neighborhoods $\neighborU$ of $\bar{x}$ and $\neighborV$ of $\bar{y}$ such that
    \begin{align*}
    \dist_S(x,F^{-1}\bar{y})\leq\kappa \dist_S(\bar{y},Fx\cap \neighborV)\qquad{\hbox{for all }} x\in\neighborU
    \end{align*}
    and some $S\in\PD$. 
  \label{def:metric_subreg}
\end{defin}
The definition is equivalent to that in \cite{Dontchev2004Regularity} but uses the $\|\cdot\|_S$ norm distance instead of the equivalent canonical norm distance to somewhat reduce notation in our convergence proof. Metric subregularity is weaker than metric regularity and the much stronger property strong monotonicity. See \cite{Dontchev2004Regularity} and \cite[Chapter~3]{Dontchev2009Implicit} for more on metric subregularity and its relation to other regularity properties. A similar linear convergence rate result under metric subregularity is provided in \cite[Theorem~3.3]{Latafat2017Asymmetric} for the AFBA algorithm that we in Section~\ref{sec:AFBA} will show is a special case of Algorithm~\ref{alg:NOFOB}. Our proof in this more general setting is similar to that for \cite[Theorem~3.3]{Latafat2017Asymmetric}.
\begin{thm}
  Suppose that Assumptions~\ref{ass:prob} and~\ref{ass:alg} hold and that the relaxation parameter $\theta_k\in[\epsilon_\theta,2-\epsilon_\theta]$ for some $\epsilon_\theta\in(0,1)$. Further suppose that $A+C$ is metrically subregular at all $z\in\zero(A+C)$ for 0 and that $\hilbert$ is either finite-dimensional or that the neighborhood $\neighborU=\hilbert$ at all $z\in\zero(A+C)$ in the metric subregularity definition. Then the $\|\cdot\|_S$-norm distance from $x_k$ to the solution set, i.e., $\dist_S(x_k,\zero(A+C))$, converges locally Q-linearly to zero and $x_k\to\bar{x}\in\zero(A+C)$ locally R-linearly.
  \label{thm:lin_conv}
\end{thm}
\begin{pf}
  We let $\neighborU_{z}$ and $\neighborV_z$ denote the neighborhoods in the metric subregularity definition (Definition~\ref{def:metric_subreg}) at the different $z\in\zero(A+C)$ for 0. Since $(A+C)^{-1}(0)=\zero(A+C)$, we conclude from Definition~\ref{def:metric_subreg} that metric subregularity of $A+C$ at all $z\in\zero(A+C)$ for 0 implies that
  \begin{align}
    \dist_S(x,\zero(A+C))\leq\kappa\dist_S(0,(A+C)x\cap\neighborV_z)\leq \kappa\|v\|_S
    \label{eq:metric_subreg_pf}
  \end{align}
  for all $x\in\neighborU_z$ and $v\in(A+C)x$ with $v\in\neighborV_z$, i.e., $\|v\|_S< \nu_z$ for some $\nu_z\in[\underline{\nu},\infty)$ where $\underline{\nu}>0$. By letting $\neighborU^\star=\cup_{z\in\zero(A+C)}\neighborU_{z}$ be a neighborhood of $\zero(A+C)$ and $\nu^\star=\inf_{z\in\zero(A+C)} \nu_z\geq\underline{\nu}>0$, we conclude that \eqref{eq:metric_subreg_pf} holds for all $x\in\neighborU^\star$ and $v\in(A+C)x$ with $\|v\|_S< \nu^\star$.

  Now, recall from \eqref{eq:inGraph} that $u_k\in(A+C)\hat{x}_k$, where $u_{k}=M_{k}x_{k}-M_{k}\hat{x}_{k}-(Cx_{k}-C\hat{x}_{k})$. Since $M_k$ is $L_M\sqrt{\tfrac{\lambda_{\max}(S^{-1})}{\lambda_{\min}(S)}}$-Lipschitz continuous  w.r.t. $\|\cdot\|_S$ (Assumption~\ref{ass:alg} and Proposition~\ref{prop:Lip}) and $C$ is $\tfrac{\beta\lambda_{\max}(S^{-1})}{\lambda_{\min}(P^{-1})}$-Lipschitz continuous w.r.t. $\|\cdot\|_S$ (Assumption~\ref{ass:prob} and Proposition~\ref{prop:coco_to_Lip}), we conclude by the triangle inequality that
   \begin{align*}
     \|u_{k}\|_S&=\|M_{k}x_{k}-M_{k}\hat{x}_{k}-(Cx_{k}-C\hat{x}_{k})\|_S\leq \xi\|x_k-\hat{x}_k\|_S\to 0
   \end{align*}
   as $k\to\infty$, where $\xi:=\left(L_M\sqrt{\tfrac{\lambda_{\max}(S^{-1})}{\lambda_{\min}(S)}}+\tfrac{\beta\lambda_{\max}(S^{-1})}{\lambda_{\min}(P^{-1})}\right)>0$ and convergence to 0 follows from Theorem~\ref{thm:conv}{\it(ii)}. If $\mathcal{U^\star}=\hilbert$, we have $\hat{x}_k\in\mathcal{U^\star}$ for all $k\in\naturals$ and if $\hilbert$ is finite-dimensional, Theorem~\ref{thm:conv}{\it(ii)} and {\it(iii)} imply that there exists a $K\in\naturals$ such that $\hat{x}_k\in\mathcal{U^\star}$ for all $k\geq K$, since weak convergence implies strong convergence in finite-dimensional settings. Therefore, \eqref{eq:metric_subreg_pf} implies that
   \begin{align*}
     \dist_S(\hat{x}_k,\zero(A+C))\leq\kappa\|u_k\|_S\leq\kappa\xi\|x_k-\hat{x}_k\|_S,
   \end{align*}
   which, using the definition of the distance function via projections and the triangle inequality, implies that
   \begin{align}
     \nonumber\dist_S(x_k,\zero(A+C))&\leq\|x_k-\Pi_{\zero(A+C)}^S\hat{x}_k\|_S\\
     \nonumber                              &\leq\|x_k-\hat{x}_k\|_S+\|\hat{x}_k-\Pi_{\zero(A+C)}^S\hat{x}_k\|_S\\
   \nonumber                       &=\|x_k-\hat{x}_k\|_S+\dist_S(\hat{x}_k,\zero(A+C))\\
                          &\leq(1+\kappa\xi)\|x_k-\hat{x}_k\|_S
                            \label{eq:dist_resid}
   \end{align}
   for all $k\geq K$. Now, let $z=\proj_{\zero(A+C)}^S(x_k)$ and use \eqref{eq:seq-Fejer} and that $\theta_k\in[\epsilon_\theta,2-\epsilon_\theta]$ followed by \eqref{eq:proj_diff_to_residual} with $\chi:=(1-\tfrac{\beta}{4})\tfrac{\lambda_{\min}(P)}{L_M^2\lambda_{\max}(S^{-1})}\tfrac{\lambda_{\min}(P)}{\lambda_{\max}(S)}>0$ and \eqref{eq:dist_resid} to conclude that
  \begin{align}
    \nonumber\dist_S^2(x_{k+1},\zero(A+C))&\leq\|x_{k+1}-z\|_S^2\\
\nonumber             &\leq\|x_k-z\|_S^2-\theta_k(2-\theta_k)\|x_k-\Pi_{H_k}^Sx_{k}\|_S^2\\
\nonumber             &\leq\|x_k-z\|_S^2-\epsilon_\theta(2-\epsilon_\theta)\|x_k-\Pi_{H_k}^Sx_{k}\|_S^2\\
\nonumber           &\leq \|x_k-z\|_S^2-\epsilon_\theta(2-\epsilon_\theta)\chi^2\|x_k-\hat{x}_k\|_{S}^2\\
\nonumber           &=\dist_S^2(x_k,\zero(A+C))-\epsilon_\theta(2-\epsilon_\theta)\chi^2\|x_k-\hat{x}_k\|_{S}^2\\
           &\leq \left(1-\tfrac{\epsilon_\theta(2-\epsilon_\theta)\chi^2}{(1+\kappa\xi)^2}\right)\dist_S^2(x_k,\zero(A+C)).\label{eq:Q_lin_dist}
  \end{align}
  We conclude Q-linear convergence since $\epsilon_\theta\in(0,1)$, $\chi>0$, $\kappa\geq 0$, and $\xi>0$.

  Let us now prove strong and R-linear convergence of $(x_k)_{k\in\naturals}$ towards the weak limit point $\bar{x}\in\zero(A+C)$ (see Theorem~\ref{thm:conv}). Let
  \begin{align*}
c:=\sqrt{\left(1-\tfrac{\epsilon_\theta(2-\epsilon_\theta)\chi^2}{(1+\kappa\xi)^2}\right)} \qquad{\hbox{and}}\qquad  R:=c^{-K}\sqrt{\tfrac{2-\epsilon_\theta}{\epsilon_\theta}}\dist_S(x_K,\zero(A+C))
  \end{align*}
  where $R\geq 0$ and $c\in[0,1)$ is the Q-linear convergence factor for $\dist_S(x_k,\zero(A+C))$ from \eqref{eq:Q_lin_dist} that holds for all $k\geq K$. Then Step~\ref{alg-step:rel_proj} in Algorithm~\ref{alg:NOFOB} and \eqref{eq:seq-Fejer} with $z=\Pi_{\zero(A+C)}^S(x_k)$ imply that for all $k\geq K$:
  \begin{align*}
    \|x_{k+1}-x_k\|_S&=\theta_k\|x_k-\Pi_{H_k}^Sx_k\|_S\leq\sqrt{\tfrac{\theta_k}{2-\theta_k}} \dist_S(x_k,\zero(A+C))\\
                     &\leq \sqrt{\tfrac{2-\epsilon_\theta}{\epsilon_\theta}} \dist_S(x_k,\zero(A+C))
                       \leq \sqrt{\tfrac{2-\epsilon_\theta}{\epsilon_\theta}} c^{k-K}\dist_S(x_K,\zero(A+C))\\
    &=Rc^{k}.
  \end{align*}
By repeated application of the triangle inequality, we get for arbitrary $j>k\geq K$ that
  \begin{align}
    \|x_k-x_j\|_S&\leq \sum_{i=0}^{j-1}\|x_{k+i+1}-x_{k+i}\|_S\leq R\sum_{i=0}^{j-1}c^{k+i}\leq Rc^k\sum_{i=0}^\infty c^{i}=\frac{R}{1-c}c^k,
                 \label{eq:resid_lin_conv}
  \end{align}
  so $(x_k)_{k\in\naturals}$ is Cauchy and converges therefore. Let $x^\star$ be the limit point, then \eqref{eq:resid_lin_conv} implies
  \begin{align}
    \|x_k-x^\star\|_S=\lim_{j\to\infty}\|x_k-x_j\|_S\leq \tfrac{R}{1-c}c^k.
    \label{eq:dist_R_rate}
  \end{align}
  Therefore $\|x_k-x^\star\|_S\to 0$ as $k\to\infty$ and since $x_k\weakto\bar{x}\in\zero(A+C)$ (Theorem~\ref{thm:conv}), we conclude from \cite[Corollary~2.52]{Bauschke2017Convex} that $x^\star=\bar{x}$ and $x_k\to\bar{x}\in\zero(A+C)$. The R-linear rate is given by \eqref{eq:dist_R_rate}.
\end{pf}

\section{A Four-Operator Splitting Method}

\label{sec:FOS}

In this section, we present a novel four-operator splitting method that is a special case of NOFOB. It solves monotone inclusion problems of the form
\begin{align}
  0\in Bx+Dx+Ex+Kx
  \label{eq:FOS_prob}
\end{align}
that satisfies the following assumption.
\begin{ass}
  Assume that $B:\hilbert\to 2^\hilbert$, that $D:\hilbert\to\hilbert$ is $L_D$-Lipschitz continuous (w.r.t. $\|\cdot\|$), $E:\hilbert\to\hilbert$ is $\tfrac{1}{\beta}$-cocoercive w.r.t. $\|\cdot\|_P$ for some $P\in\PD$ and $\beta\in[0,4)$, $K:\hilbert\to\hilbert$ is linear skew-adjoint, $B+D$ is maximally monotone, and that the solution set $\zero(B+D+E+K)$ is nonempty.
  \label{ass:FOS_prob}
  \end{ass}
We only require the sum $B+D$ to be maximally monotone, not the individual operators $B$ and $D$. Since $K$ is maximally monotone \cite[Example~20.35]{Bauschke2017Convex} and has full domain, also $B+D+K$ is maximally monotone \cite[Corollary~25.5]{Bauschke2017Convex}.

We construct the four-operator splitting method from Algorithm~\ref{alg:NOFOB_explicit}. The single-valued operators $D$ and $K$ are put in $A$ and removed from the inversion by subtracting them using the kernel $M_k$. We let $C=E$ and $A=B+D+K$ and use $M_k=Q_k-D-K$ where $Q_k:\hilbert\to\hilbert$ is strongly monotone and Lipschitz continuous. The nonlinear forward-backward step in Algorithm~\ref{alg:NOFOB_explicit} becomes
\begin{align*}
  \hat{x}_k&=(M_k+A)^{-1}(M_k-C)x_k\\
           &= (Q_k-D-K+B+D+K)^{-1}(Q_k-D-K-E)x_k\\
           &=(Q_k+B)^{-1}(Q_k-D-K-E)x_k.
\end{align*}
The $D$ and $K$ operators are canceled by $M_k$ from the backward part and effectively moved to the forward part. We use this construction in Algorithm~\ref{alg:NOFOB_explicit} to arrive at the four-operator splitting method below. We see that $M_k$ must be evaluated at $\hat{x}_k$ for the projection correction step. This implies that $D$ and $K$ are evaluated twice per iteration, while $E$ and the resolvent of $B$ are evaluated once. 
\ifx\one\showalgs\begin{algorithm}[H]
  \begin{algorithmic}[1]
    \STATE {\bf{Let:}} $M_k=Q_k-D-K$
    \STATE {\bf{Input:}} $x_0\in\hilbert$
		\FOR{$k=0,1,\ldots$}
                \STATE $\hat{x}_{k}:=(Q_k+B)^{-1}(Q_k-D-K-E)x_k$\label{item:FOS_res_update}
                \STATE $\mu_k:=\frac{\langle M_kx_k-M_k\hat{x}_k,x_k-\hat{x}_k\rangle-\tfrac{\beta}{4}\|x_k-\hat{x}_k\|_P^2}{\|M_kx_k-M_k\hat{x}_k\|_{ S^{-1}}^2}$
                \STATE $x_{k+1}:=x_k-\theta_k\mu_kS^{-1}(M_kx_k-M_k\hat{x}_k)$
		\ENDFOR
	\end{algorithmic}
        \caption{Four-Operator Splitting}
        \label{alg:FOS}
      \end{algorithm}\fi

      Since the algorithm is a direct special case of Algorithm~\ref{alg:NOFOB_explicit}, we get convergence from Corollary~\ref{cor:conv_exp} (and linear convergence from Theorem~\ref{thm:lin_conv} under metric subregularity) if $Q_k$ and $P\in\PD$ are chosen such that Assumptions~\ref{ass:prob}~and~\ref{ass:alg} are satisfied.
      \begin{cor}
        Suppose that Assumption~\ref{ass:FOS_prob} holds, that for each $k\in\naturals$ we have $\theta_k\in[\epsilon_\theta,2-\epsilon_\theta]$ for some $\epsilon_\theta\in(0,1)$, and that $Q_k:\hilbert\to\hilbert$ is maximally monotone, single-valued, and $L_Q$-Lipschitz continuous for some $L_Q>0$, and that $Q_k-D$ is 1-strongly monotone w.r.t. $\|\cdot\|_P$. Then Algorithm~\ref{alg:FOS} generates a sequence $(x_k)_{k\in\naturals}$ such that $x_k\weakto \bar{x}\in\zero(B+D+E+K)$. Moreover, if $B+D+E+K$ is metrically subregular (Definition~\ref{def:metric_subreg}) at all $z\in\zero(B+D+E+K)$ for 0 and $\hilbert$ is finite-dimensional, then $x_k\to\bar{x}\in\zero(B+D+E+K)$ locally R-linearly and $\dist_S(x_k,\zero(B+D+E+K))$ converges locally Q-linearly to zero.
        \label{cor:FOS_conv}
      \end{cor}
      \begin{pf}
        Let $A=B+D+K$ and $C=E$. Assumption~\ref{ass:prob} follows from Assumption~\ref{ass:FOS_prob} by noting that $A$ is maximally monotone by \cite[Corollary~25.5]{Bauschke2017Convex} since $B+D$ and $K$ are maximally monotone \cite[Example~20.35]{Bauschke2017Convex} and $K$ has full domain. Single-valuedness and Lipschitz continuity of $M_k=Q_k-D-K$ are immediate. Strong monotonicity with modulus 1 w.r.t. $\|\cdot\|_P$ of $Q_k-D$ is equivalent to the same property for $M_k=Q_k-D-K$ since adding the linear skew-adjoint operator $-K$ does not affect strong monotonicity. The last property for $M_k$ in Assumption~\ref{ass:alg} is maximal monotonicity. It holds by \cite[Corollary~20.28]{Bauschke2017Convex} since all $M_k:\hilbert\to\hilbert$ are (strongly) monotone and (Lipschitz) continuous. Since Assumptions~\ref{ass:prob} and~\ref{ass:alg} are satisfied and Algorithm~\ref{alg:FOS} is a special case of Algorithm~\ref{alg:NOFOB_explicit}, Corollary~\ref{cor:conv_exp} and Theorem~\ref{thm:lin_conv} provide the convergence results.
      \end{pf}

      \begin{rem}
        The Lipschitz continuous operator $D$ restricts the choice of $Q_k$ via the required 1-strong monotonicity of $Q_k-D$ w.r.t. $\|\cdot\|_P$. The cocoercive operator $E$ also restricts the choice since it is cocoercive w.r.t. the same norm $\|\cdot\|_P$ and $P\in\PD$ must be selected to satisfy both conditions. However, the linear skew-adjoint operator $K$ does not restrict the choice of $Q_k$. If, for instance, $D=E=0$, $P\in\PD$ is arbitrary and can be chosen $P=\epsilon\id$ for arbitrary $\epsilon\in(0,1)$ and the only conditions we get on all $Q_k$ are $\epsilon$-strong monotonicity and $\tfrac{1}{\epsilon}$-Lipschitz continuity. 
      \end{rem}

      In the remainder of this section, we will state special cases of the four-operator splitting method. We will first consider cases with $Q_k=\gamma_k^{-1}\id$ that gives forward-backward-forward splitting (if $K=E=0$) and forward-backward-half-forward splitting (if $K=0$) as special cases. We then let $Q_k=P+G$, where $G$ is linear skew-adjoint, to get asymmetric forward-backward-adjoint splitting as a special case (if $D=0$). We will finish the section with a discussion on standard forward-backward splitting (with $K=D=0$) and comment on consequences of nonsymmetry or nonlinearity of $M_k$, or lack thereof, on the projection correction step.

      \subsection{Special cases with $Q_k=\gamma_k^{-1}\id$}

      We will let $Q_k=\gamma_k^{-1}\id$ with $\gamma_k>0$ and derive explicit conditions on $\gamma_k$ that guarantee convergence. This choice gives $M_k=\gamma_k^{-1}\id-D-K$ and nonlinear forward-backward update
      \begin{align}
        \nonumber\hat{x}_k&=(\gamma_k^{-1}\id+B)^{-1}(\gamma_k^{-1}\id-D-K-E)x_k\\
        &=(\id+\gamma_kB)^{-1}(\id-\gamma_k(D+K+E))x_k\label{eq:FOS_ID_FB_update}
      \end{align}
      in Step~\ref{item:FOS_res_update} of Algorithm~\ref{alg:FOS}. We see that $\gamma_k$ acts as a step-size in the update. The resulting special case of Algorithm~\ref{alg:FOS} is stated below.
\ifx\one\showalgs\begin{algorithm}[H]
  \begin{algorithmic}[1]
    \STATE {\bf{Let:}} $M_k=\gamma_k^{-1}\id-D-K$
    \STATE {\bf{Input:}} $x_0\in\hilbert$
		\FOR{$k=0,1,\ldots$}
                \STATE $\hat{x}_{k}:=(\id+\gamma_kB)^{-1}(\id-\gamma_k(D+K+E))x_k$\label{item:FOS_res_update_gammak}
                \STATE $\mu_k:=\frac{\langle M_kx_k-M_k\hat{x}_k,x_k-\hat{x}_k\rangle-\tfrac{\beta_E}{4}\|x_k-\hat{x}_k\|^2}{\|M_kx_k-M_k\hat{x}_k\|_{ S^{-1}}^2}$\label{item:muk_update}
                \STATE $x_{k+1}:=x_k-\theta_k\mu_kS^{-1}(M_kx_k-M_k\hat{x}_k)$
		\ENDFOR
	\end{algorithmic}
        \caption{Four-Operator Splitting with $Q_k=\gamma_k^{-1}\id$}
        \label{alg:FOS_gammak}
      \end{algorithm}\fi

   We will also analyze a short-step version that has standard formulations of forward-backward-forward and forward-backward-half-forward as special cases. 
      We let $S=\id$ and replace $\mu_k$ by $\hat{\mu}_k$ (that satisfies \eqref{eq:mu_bound} for all $x\neq y\in\hilbert$ and all $k\in\naturals$) and set $\theta_k=\tfrac{\gamma_k}{\hat{\mu}_k}$ in Algorithm~\ref{alg:FOS_gammak} to get the following algorithm.
\ifx\one\showalgs\begin{algorithm}[H]
  \begin{algorithmic}[1]
\STATE {\bf Input:} $x_0\in\hilbert$
    \FOR{$k=0,1,\ldots$}
                \STATE $\hat{x}_{k}:=(\id+\gamma_k B)^{-1}(\id-\gamma_k(D+K+E))x_k$ 
                \STATE $x_{k+1}:=\hat{x}_k-\gamma_k((D+K)\hat{x}_k- (D+K)x_k)$
		\ENDFOR
	\end{algorithmic}
        \caption{Four-Operator Splitting with $Q_k=\gamma_k^{-1}\id$ and Conservative Step-Length}
        \label{alg:FOS_conservative}
      \end{algorithm}\fi

      Algorithm~\ref{alg:FOS_conservative} is a conservative version of Algorithm~\ref{alg:FOS_gammak} since in Algorithm~\ref{alg:FOS_conservative}, the step-length is decided by the global constant $\hat{\mu}_k$ that satisfies \eqref{eq:mu_bound} for all $x\neq y\in\hilbert$, while in Algorithm~\ref{alg:FOS_gammak}, the step-length is decided by a local constant $\mu_k$ that satisfies \eqref{eq:muk}, which is \eqref{eq:mu_bound} but only for $x_k$ and $\hat{x}_k$. Algorithm~\ref{alg:FOS_conservative} also uses a specific relaxation parameter $\theta_k=\tfrac{\gamma_k}{\hat{\mu}_k}$, while Algorithm~\ref{alg:FOS_gammak} can use arbitrary $\theta\in[\epsilon_\theta,2-\epsilon_\theta]$. The requirement for convergence in Corollary~\ref{cor:FOS_conv} that $\theta_k\in[\epsilon_\theta,2-\epsilon_\theta]$ sets additional requirements on the step-sizes $\gamma_k$ in Algorithm~\ref{alg:FOS_conservative} due to the choice $\theta_k=\tfrac{\gamma_k}{\hat{\mu}_k}$. We will in the convergence analysis in Section~\ref{sec:FOS_gammak_conv} see that the allowed range of step-sizes $\gamma_k$ is larger for Algorithm~\ref{alg:FOS_gammak} than for Algorithm~\ref{alg:FOS_conservative}, and that it does not depend on $K$ in Algorithm~\ref{alg:FOS_gammak}, while it does in Algorithm~\ref{alg:FOS_conservative}.

\subsubsection{Convergence Analysis}
\label{sec:FOS_gammak_conv}

We will in this section analyze Algorithms~\ref{alg:FOS_gammak} and \ref{alg:FOS_conservative}. We will give conditions on $\gamma_k$ for convergence under the following assumption that differs only in the assumption on the cocoercive term $E$ compared to Assumption~\ref{ass:FOS_prob}.
\begin{ass}
        Assume that $B:\hilbert\to 2^\hilbert$, that $D:\hilbert\to\hilbert$ is $L_D$-Lipschitz continuous (w.r.t. $\|\cdot\|$), $E:\hilbert\to\hilbert$ is $\tfrac{1}{\beta_E}$-cocoercive w.r.t. $\|\cdot\|$, $K:\hilbert\to\hilbert$ is linear skew-adjoint, $B+D$ is maximally monotone, and that the solution set $\zero(B+D+E+K)$ is nonempty.
  \label{ass:FOS_prob_gammak}
\end{ass}

Before we proceed with the convergence analysis, let us define
\begin{align}
  P:=(\underline{\gamma}^{-1}-L_D)\id\qquad{\hbox{where}}\qquad
  \underline{\gamma}^{-1}:=\inf_{k\in\naturals}\gamma_k^{-1},
  \label{eq:Pdef_gamma_min}
\end{align}
and note that Step~\ref{item:muk_update} in Algorithm~\ref{alg:FOS_gammak} has been equivalently been rewritten compared to the same step in Algorithm~\ref{alg:FOS} with the equality $\tfrac{\beta_E}{4}\|x_k-\hat{x}_k\|^2=\tfrac{\beta}{4}\|x_k-\hat{x}_k\|_P^2$ that holds under Assumption~\ref{ass:FOS_prob_gammak} due to the following result.
\begin{lem}
          Assume that $E:\hilbert\to\hilbert$ is $\tfrac{1}{\beta_E}$-cocoercive w.r.t. $\|\cdot\|$ and let $P=(\underline{\gamma}^{-1}-L_D)\id\in\PD$. Then $E$ is $\tfrac{1}{\beta}$-cocoercive w.r.t. $\|\cdot\|_P$ with $\beta=\tfrac{\beta_E}{\underline{\gamma}^{-1}-L_D}$. 
        \label{lem:C_coco_FOS}
      \end{lem}
      \begin{pf}
        Let $\delta:=\underline{\gamma}^{-1}-L_D$ that satisfies $\delta>0$ since $P=\delta\id\in\PD$. Since $E$ is $\tfrac{1}{\beta_E}$ cocoercive w.r.t $\|\cdot\|$, we have
      \begin{align*}
        \langle Cx-Cy,x-y\rangle&\geq\tfrac{1}{\beta_E}\|Cx-Cy\|^2=\tfrac{\delta}{\beta_E}\|Cx-Cy\|_{\delta^{-1}\id}^2
                                =\tfrac{\delta}{\beta_E}\|Cx-Cy\|_{P^{-1}}^2
      \end{align*}
      for all $x,y\in\hilbert$.
    \end{pf}

      Let us also state the following lemma that will help us verify Assumption~\ref{ass:alg}. 
      The proof of the lemma is deferred to Appendix~\ref{app:M_prop_FOS_pf}.

\begin{lem}
  Suppose that Assumption~\ref{ass:FOS_prob_gammak} holds and let $P=(\underline{\gamma}^{-1}-L_D)\id\in\PD$. Then all $M_k=\gamma_k^{-1}\id-D-K$ with $k\in\naturals$ are single-valued, maximally monotone, 1-strongly monotone w.r.t. $\|\cdot\|_P$, and $(\gamma_k^{-1}+L_D+\|K\|)$- Lipschitz continuous w.r.t. $\|\cdot\|$. 
\label{lem:M_prop_FOS}
\end{lem}      
We are ready to show convergence of the four operator splitting method in Algorithm~\ref{alg:FOS_gammak}.
          \begin{thm}
            Suppose that Assumption~\ref{ass:FOS_prob_gammak} holds, that $\theta_k\in[\epsilon_\theta,2-\epsilon_\theta]$ for some $\epsilon_\theta\in(0,1)$, and that $\gamma_{k}\in[\epsilon,\tfrac{1}{\epsilon}]$ and 
            $\gamma_k\leq\tfrac{4-\epsilon}{\beta_E+4L_D}$ for some $\epsilon\in(0,1)$. Then Algorithm~\ref{alg:FOS_gammak} generates a sequence $(x_k)_{k\in\naturals}$ such that $x_k\weakto \bar{x}\in\zero(B+D+E+K)$. Moreover, if $B+D+E+K$ is metrically subregular (Definition~\ref{def:metric_subreg}) at all $z\in\zero(B+D+E+K)$ for 0 and $\hilbert$ is finite-dimensional, then $x_k\to\bar{x}\in\zero(B+D+E+K)$ locally R-linearly and $\dist_S(x_k,\zero(B+D+E+K))$ converges locally Q-linearly to zero.
      \label{thm:FOS_conv_gammak}
    \end{thm}
    \begin{pf}
      Since Algorithm~\ref{alg:FOS_gammak} is a special case of Algorithm~\ref{alg:NOFOB_explicit} with $A=B+D+K$, $C=E$, and $M_k=\gamma_k^{-1}\id-D-K$, Corollary~\ref{cor:conv_exp} and Theorem~\ref{thm:lin_conv} give the convergence results if Assumptions~\ref{ass:prob} and~\ref{ass:alg} hold. We let $P=(\underline{\gamma}^{-1}-L_D)\id$ and note that $P\in\PD$ if $\underline{\gamma}^{-1}=\inf_{k\in\naturals}\gamma_k^{-1}>L_D$. For $L_D>0$ it holds since $\gamma_k\leq\tfrac{4-\epsilon}{\beta_E+4L_D}\leq\tfrac{1-\epsilon/4}{L_D}$ and for $L_D=0$ it holds since $\gamma_k\leq\tfrac{1}{\epsilon}$. 

      For Assumption~\ref{ass:prob}, we only need to verify {\it{(\ref{ass:prob:C})}} since {\it{(\ref{ass:prob:A})}} is proven in the proof of Corollary~\ref{cor:FOS_conv} and {\it{(\ref{ass:prob:nonempty})}} trivially follows from Assumption~\ref{ass:FOS_prob_gammak}.
      
      To verify Assumption~\ref{ass:prob}{\it{(\ref{ass:prob:C})}}, we note that Lemma~\ref{lem:C_coco_FOS} implies that $C=E$ is $\tfrac{1}{\beta}$-cocoercive w.r.t. $\|\cdot\|_P$ with  $\beta=\tfrac{\beta_E}{\underline{\gamma}^{-1}-L_D}$. The denominator in the cocoercivity constant $\tfrac{\beta_E}{\underline{\gamma}^{-1}-L_D}$  is positive since $0<\gamma_k\leq\tfrac{1-\epsilon/4}{L_D}$, so it remains to show that $\tfrac{\beta_E}{\underline{\gamma}^{-1}-L_D}<4$. Since $\underline{\gamma}^{-1}-L_D>0$ and $\beta_E>0$, we have
      \begin{align*}
        \tfrac{\beta_E}{\underline{\gamma}^{-1}-L_D}=\tfrac{\beta_E}{\inf_{k\in\naturals}\gamma_k^{-1}-L_D}=\sup_{k\in\naturals}\tfrac{\beta_E}{\gamma_k^{-1}-L_D}\leq\tfrac{\beta_E}{\tfrac{\beta_E+4L_D}{4-\epsilon}-L_D}=\tfrac{\beta_E(4-\epsilon)}{\beta_E+L_D\epsilon} \leq 4-\epsilon,
      \end{align*}
where we have used $\gamma_k\leq\tfrac{4-\epsilon}{\beta_E+4L_D}$ in the first inequality.


That Assumption~\ref{ass:alg} holds follows from Lemma~\ref{lem:M_prop_FOS} by defining $L_M=\epsilon^{-1}+L_D+\|K\|$ for Lipschitz continuity and noting that $\gamma_k^{-1}+L_D+\|K\|\leq\epsilon^{-1}+L_D+\|K\|$ for all $k\in\naturals$ since $\gamma_k\geq\epsilon>0$.
    \end{pf}
    \begin{rem}
      Note that the linear skew-adjoint operator $K$ does not restrict the choice of step-size $\gamma_k$. In particular, if $D=E=0$, the only step-size requirement in Theorem~\ref{thm:FOS_conv_gammak} is that $\gamma_k\in[\epsilon,\tfrac{1}{\epsilon}]$ for arbitrary $\epsilon\in(0,1)$. We will in Section~\ref{sec:Projective} see that synchronous projective splitting is a special case of this setting. 
    \end{rem}

      To prove convergence of Algorithm~\ref{alg:FOS_conservative} via Corollary~\ref{cor:conv_conserv}, we; provide a $\hat{\mu}_k$ for all $k\in\naturals$ that satisfies $\hat{\mu}_k\geq\epsilon_{\mu}$ for some $\epsilon_{\mu}>0$ and satisfies \eqref{eq:mu_bound} with $S=\id$ for all $x\neq y\in\hilbert$, and show that the relaxation parameters $\theta_k=\tfrac{\gamma_k}{\hat{\mu}_k}$ satisfy the conditions $\theta_k\in(0,2)$ and $\liminf_{k\to\infty}\theta_k(2-\theta_k)>0$. We will use the following proposition, that is proven in Appendix~\ref{app:mu_bound_FOS}, to show this.
    \begin{prop}
      Suppose that Assumption~\ref{ass:FOS_prob_gammak} holds and that the step-size $\gamma$ satisfies $\gamma\in[\epsilon,\tfrac{1}{\epsilon}]$ and $\gamma\leq\tfrac{4-\epsilon}{\beta_E+\sqrt{\beta_E^2+16(L_D+\|K\|)^2}}$ for some $\epsilon\in(0,1)$ that satisfies $\epsilon^{-1}\geq\tfrac{\beta_EL_D\epsilon}{2(1-\epsilon)}$ (which trivially holds for small enough $\epsilon>0$). Let
\begin{align}
  \bar{\epsilon}&=\frac{\epsilon((8-\epsilon)\sqrt{\beta_E^2+16(L_D+\|K\|)^2}+\epsilon\beta_E)}{4(8-\epsilon)}, \label{eq:ep_bar_def}\\
  \delta&=\frac{\bar{\epsilon}}{2(\epsilon^{-1}+L_D+\|K\|-\tfrac{\beta_EL_D\epsilon}{4(1-\epsilon)})}.
          \label{eq:delta_def}
\end{align}
    Then $\delta\in(0,1)$ and $M=\gamma^{-1}\id-D-K$ satisfies
      \begin{align}
        \frac{\gamma}{2-\delta}\leq\frac{\langle Mx-My,x-y\rangle-\tfrac{\beta}{4}\|x-y\|_P^2}{\|Mx-My\|^2}
        \label{eq:gamma_over_twodelta}
      \end{align}
      with $P=(\gamma^{-1}-L_D)\id$ for all $x,y\in\hilbert$ such that $x\neq y$.
      \label{prop:mu_bound_FOS}
    \end{prop}
We are ready to state our convergence result for Algorithm~\ref{alg:FOS_conservative}.
\begin{thm}
  Suppose that Assumption~\ref{ass:FOS_prob_gammak} holds, that $\theta_k\in[\epsilon_\theta,2-\epsilon_\theta]$ for some $\epsilon_\theta\in(0,1)$, and that $\gamma_{k}\in[\epsilon,\tfrac{1}{\epsilon}]$ and 
            $\gamma_k\leq\tfrac{4-\epsilon}{\beta_E+\sqrt{\beta_E^2+16(L_D+\|K\|)^2}}$ for some $\epsilon\in(0,1)$. Then Algorithm~\ref{alg:FOS_conservative} generates a sequence $(x_k)_{k\in\naturals}$ such that $x_k\weakto \bar{x}\in\zero(B+D+E+K)$. Moreover, if $B+D+E+K$ is metrically subregular (Definition~\ref{def:metric_subreg}) at all $z\in\zero(B+D+E+K)$ for 0 and $\hilbert$ is finite-dimensional, then $x_k\to\bar{x}\in\zero(B+D+E+K)$ locally R-linearly and $\dist_S(x_k,\zero(B+D+E+K))$ converges locally Q-linearly to zero.
      \label{thm:FOS_conv_conservative}
    \end{thm}
    \begin{pf}
      Algorithm~\ref{alg:FOS_conservative} is a special case of Algorithm~\ref{alg:NOFOB_constant_step} with $A=B+D+K$, $C=E$, $M_k=\gamma_k^{-1}\id-D-K$, and $\theta_k=\tfrac{\gamma_k}{\hat{\mu}_k}$. Therefore, Corollary~\ref{cor:conv_conserv} and Theorem~\ref{thm:lin_conv}, give the convergence results if Assumptions~\ref{ass:prob} and~\ref{ass:alg} hold, if there exists $\epsilon_{\mu}>0$ and $\hat{\mu}_k$ such that for all $k\in\naturals$, $\hat{\mu}_k\geq\epsilon_{\mu}$ and $\hat{\mu}_k$ satisfies \eqref{eq:mu_bound} for all $x,y\in\hilbert$ with $x\neq y$, and if $\theta_k=\tfrac{\gamma_k}{\hat{\mu}_k}\in[\epsilon_\theta,2-\epsilon_\theta]$ for some $\epsilon_\theta\in(0,1)$. 
      
      That Assumptions~\ref{ass:prob} and \ref{ass:alg} hold is proven in the same way as for Theorem~\ref{thm:FOS_conv_gammak} since the upper bound $\gamma_k\leq\tfrac{4-\epsilon}{\beta_E+\sqrt{\beta_E^2+16(L_D+\|K\|)^2}}$ is smaller than the upper bound in Theorem~\ref{thm:FOS_conv_gammak}, $\gamma_k\leq\tfrac{4-\epsilon}{\beta_E+4L_D}$. 

      We let $\hat{\mu}_k=\tfrac{\gamma_k}{2-\delta}$, where $\delta$ is defined in \eqref{eq:delta_def}. Since $\gamma_k\geq \epsilon$ and $\delta\in(0,1)$ is iteration independent (Proposition~\ref{prop:mu_bound_FOS}), we conclude that $\hat{\mu}_k\geq\epsilon_{\mu}:=\tfrac{\epsilon}{2-\delta}$. Proposition~\ref{prop:mu_bound_FOS} also implies that $\hat{\mu}_k$ satisfies \eqref{eq:mu_bound} for all $x,y\in\hilbert$ such that $x\neq y$. 
      
      Finally $\theta_k=\tfrac{\gamma_k}{\hat{\mu}_k}=2-\delta\in[\epsilon_\theta,2-\epsilon_\theta]$ with $\epsilon_\theta=\delta\in(0,1)$ since $\delta$ is iteration independent (Proposition~\ref{prop:mu_bound_FOS}).
    \end{pf}



    \subsubsection{Forward-Backward-Forward Splitting}

    \label{sec:FBF}

The forward-backward-forward splitting algorithm (FBF) in \cite{Tseng2000AModified} solves monotone inclusion problems of the form
\begin{align}
  0\in Bx+Dx,
  \label{eq:FBF_prob}
\end{align}
where $B:\hilbert\to 2^\hilbert$ and $D:\hilbert\to\hilbert$ are maximally monotone and $D$ is $L_D$-Lipschitz continuous. Forward-backward-forward splitting can be recovered from Algorithm~\ref{alg:FOS} by letting $E=K=0$. That the cocoercive term $E=0$ implies that it is a special case of the nonlinear backward (resolvent) method and that the sum $B+D$ is treated as one operator. Therefore, we require only the sum $A=B+D$ to be maximally monotone, not the individual operators.

The standard FBF formulation in \cite{Tseng2000AModified} is obtained by letting $E=K=0$ in Algorithm~\ref{alg:FOS_conservative} and the standard (albeit rewritten) step-size restriction $\gamma_k\in[\epsilon,\tfrac{1-\epsilon/4}{L_D}]$ under Lipschitz continuity of $D$ follows from Theorem~\ref{thm:FOS_conv_conservative}. This is a conservative version of Algorithm~\ref{alg:FOS_gammak}. A long-step FBF variation based on Algorithm~\ref{alg:FOS_gammak} with $K=E=0$ can also be found in \cite{Tseng2000AModified} by utilizing the iteration-dependent set that contains the solution set. Convergence and and linear convergence under metric subregularity of $B+D$ follow from both these methods by Theorem~\ref{thm:FOS_conv_conservative} and Theorem~\ref{thm:FOS_conv_gammak}.

\subsubsection{Forward-Backward-Half-Forward Splitting}
\label{sec:FBHF}

The forward-backward-half-forward splitting method (FBHF) in \cite{BricenoArias2018Forward} solves monotone inclusion problems of the form
\begin{align*}
  0\in Bx+Dx+Ex,
\end{align*}
where $B:\hilbert\to 2^\hilbert$, $D:\hilbert\to\hilbert$ is $L_D$-Lipschitz continuous (w.r.t. $\|\cdot\|$), $E:\hilbert\to\hilbert$ is $\tfrac{1}{\beta_E}$-cocoercive w.r.t. $\|\cdot\|$, and $B+D$ is maximally monotone. This adds a cocoercive term $E$ compared to FBF and is a special case of \eqref{eq:FOS_prob} with $K=0$.

The FBHF formulation treated in \cite[Theorem~2.3]{BricenoArias2018Forward} is exactly Algorithm~\ref{alg:FOS_conservative} with $K=0$. Theorem~\ref{thm:FOS_conv_conservative} uses the same step-size conditions for convergence as \cite[Theorem~2.3]{BricenoArias2018Forward} and shows linear convergence under metric subregularity of $B+D+E$, which appears to be a new result. A novel long-step variation of FBHF is obtained by letting $K=0$ in Algorithm~\ref{alg:FOS_gammak}. We show in Theorem~\ref{thm:FOS_conv_gammak} convergence under less restrictive step-size assumptions than in Theorem~\ref{thm:FOS_conv_conservative} and \cite[Theorem~2.3]{BricenoArias2018Forward} and show linear convergence under metric subregularity of $B+D+E$.

\subsection{Asymmetric Forward-Backward-Adjoint Splitting}

\label{sec:AFBA}
The AFBA algorithm -- asymmetric forward-backward-adjoint splitting -- in \cite{Latafat2017Asymmetric} solves monotone inclusion problems of the form
\begin{align*}
0\in Bx+Ex+Kx,
\end{align*}
where $B:\hilbert\to 2^\hilbert$ is maximally monotone, $E:\hilbert\to\hilbert$ is $\tfrac{1}{\beta}$-cocoercive w.r.t. $\|\cdot\|_P$ for some $P\in\PD$ and $\beta\in[0,4)$, and $K:\hilbert\to\hilbert$ is linear skew-adjoint.

We arrive at the AFBA algorithm in \cite{Latafat2017Asymmetric} from Algorithm~\ref{alg:FOS} by letting $D=0$ and $Q_k=Q=P+G$, where $P\in\PD$ and $G:\hilbert\to\hilbert$ is linear skew-adjoint. With this choice of $Q_k$, it is straight-forward to verify the assumptions in Corollary~\ref{cor:FOS_conv} that recovers the convergence and linear convergence rate results in \cite[Theorems~3.1~and~3.3]{Latafat2017Asymmetric}.

A special case of the fixed step-size version of AFBA in \cite{Latafat2018PrimalDual} is recovered by letting $\theta_k\mu_k=1$. We recover the corresponding convergence results in \cite{Latafat2018PrimalDual} from Corollary~\ref{cor:FOS_conv} if $\theta_k=\tfrac{1}{\mu_k}\in[\epsilon_\theta,2-\epsilon_\theta]$ for some $\epsilon_\theta\in(0,1)$. This requires that $\mu_k$ is bounded (Proposition~\ref{prop:muk_bounds}) and satisfies $\mu_k\geq\tfrac{1}{2-\epsilon_\theta}$. This is true if $(1-\tfrac{\beta}{4})P-\tfrac{1}{2-\epsilon_\theta}(Q-K^*)S^{-1}(Q-K)\succeq 0$ (which is the same condition as in \cite{Latafat2018PrimalDual} in this setting) since then
\begin{align*}
  \mu_k&=\frac{\langle(Q-K)x_k-(Q-K)\hat{x}_k,x_k-\hat{x}_k\rangle-\tfrac{\beta}{4}\|x_k-\hat{x}_k\|_P^2}{\|(Q-K)x_k-(Q-K)\hat{x}_k\|_{S^{-1}}^2}\\
       &=(1-\tfrac{\beta}{4})\frac{\|x_k-\hat{x}_k\|_P^2}{\|x_k-\hat{x}_k\|^2}_{(Q-K^*)S^{-1}(Q-K)}\geq\tfrac{1}{2-\epsilon_\theta}.
\end{align*}

\subsubsection{Special cases}

We benefit from the lists of special cases of AFBA in \cite{Latafat2017Asymmetric,Latafat2018PrimalDual} to conclude that Solodov and Tseng~\cite{Solodov1996Modified}, Douglas-Rachford~\cite{Lions1979Splitting}, ADMM \cite{Gabay1976ADual,Glowinski1975Sur}, and \cite{Drori2015ASimple,BricenoArias2011AMonotone,He2012Convergence,Latafat2019ANew}, are special cases of this algorithm (although the analysis of Douglas-Rachford and ADMM needs a slight generalization of our results with positive semidefinite $Q$).

    \subsection{Forward-Backward Splitting}
    \label{sec:FBS}
    
    Forward-backward splitting (and therefore Chambolle-Pock~\cite{Chambolle2011AFirstOrder} and V{\~u}-Condat~\cite{Condat2013APrimalDual,Vu2013ASplitting}) is a special case of Algorithm~\ref{alg:FOS} with $D=K=0$ that solves
    \begin{align*}
      0\in Bx+Ex,
    \end{align*}
    where $B:\hilbert\to 2^\hilbert$ is maximally monotone and $E:\hilbert\to\hilbert$ is $\tfrac{1}{\beta_E}$-cocoercive w.r.t. $\|\cdot\|_M$ for some $M\in\PD$. We let all $M_k=Q_k=\gamma_k^{-1}M$ and let the projection metric $S=M$. Then $S^{-1}M_k=\gamma_k^{-1}\id$ and $\mu_k=\gamma_k(1-\tfrac{\beta_E\gamma_k}{4})$ and Algorithm~\ref{alg:FOS} reduces to relaxed forward-backward splitting
    \begin{align}
      x_{k+1}&=(1-\theta_k(1-\tfrac{\beta_E\gamma_k}{4}))x_k+\theta_k(1-\tfrac{\beta_E\gamma_k}{4})(M+\gamma_kB)^{-1}(M-\gamma_kE)x_{k}.
               \label{eq:FBS_relaxed}
    \end{align}
    We let $P=\inf_{k\in\naturals}\gamma_k^{-1}M$, that satisfies $P\in\PD$ if $\gamma_k\in[\epsilon,\tfrac{1}{\epsilon}]$ for some $\epsilon\in(0,1)$. Therefore, all $M_k$ are 1-strongly monotone w.r.t. $\|\cdot\|_P$ and $E$ is $\tfrac{1}{\beta}=\tfrac{\inf_{k\in\naturals}\gamma_k^{-1}}{\beta_E}$-cocoercive w.r.t. $\|\cdot\|_P$. To satisfy Assumption~\ref{ass:prob} with $A=B$ and $C=E$, we need $\beta\in[0,4)$ which is implied by $\gamma_k\leq\tfrac{4-\epsilon}{\beta_E}$ for all $k\in\naturals$ and some $\epsilon>0$. Therefore, Corollary~\ref{cor:FOS_conv} gives convergence of \eqref{eq:FBS_relaxed} whenever, e.g., $\theta_k\in[\epsilon_\theta,2-\epsilon_\theta]$, $\gamma_k\in(\epsilon,\tfrac{1}{\epsilon})$ and $\gamma_k\leq\tfrac{4-\epsilon}{\beta_E}$. This doubles the allowed step-size range compared to standard analyses based on averaged operators \cite[Theorem~26.14]{Bauschke2017Convex}. The extended range of allowed step-size parameters has previously been shown also in \cite{Latafat2017Asymmetric} and in \cite{Giselsson2019OnCompositions} in the context of conically averaged operators.

    If $\theta_k=\tfrac{4}{4-\beta_E\gamma_k}$, then \eqref{eq:FBS_relaxed} reduces to standard (non-relaxed) forward-backward splitting
    \begin{align}
      x_{k+1}=(M+\gamma_kB)^{-1}(M-\gamma_kE)x_k.
      \label{eq:FBS_iter}
    \end{align}
    The requirement that $\theta_k\in[\epsilon_\theta,2-\epsilon_\theta]$ holds if $\gamma_k\leq\tfrac{2-\epsilon}{\beta_E}$ for some $\epsilon>0$, which is the standard requirement on $\gamma_k$. For larger $\gamma_k\geq\tfrac{2}{\beta_E}$ in the extended range, $\theta_k$ must be chosen smaller than $\tfrac{4}{4-\beta_E\gamma_k}$ to not violate $\theta_k\in[\epsilon,2-\epsilon]$. This necessitates an under-relaxed iteration \eqref{eq:FBS_relaxed} to have convergence. 

    \begin{rem}
We have selected $M_k$ to be a scalar times a positive definite self-adjoint operator $M$ and projection kernel $S=M$. As is evident from the resulting iterations \eqref{eq:FBS_relaxed} and \eqref{eq:FBS_iter}, this choice avoids the second application of $M_k$ at $\hat{x}_k$ that is used for the projection correction in Algorithm~\ref{alg:FOS}. The reason is that the projection point onto the separating hyperplane is exactly $\hat{x}_k$, which is already known. This distinguishes such choices of self-adjoint positive definite $M_k$ from nonlinear or nonsymmetric ones for which the projection correction, and hence the second $M_k$ application, is needed in every iteration to guarantee convergence. An interpretation is that our algorithm automatically adapts to nonlinearity or nonsymmetry of $M_k$. If $M_k$ is nonlinear or nonsymmetric, it makes use of the projection correction step, that needs a second evaluation of $M_k$, to guarantee convergence, while if $M$ is linear, self-adjoint, and positive definite the correction step is redundant and the second application of $M_k$ is avoided.
\end{rem}

    \ifx\showps\one
    \section{Projective Splitting}
    \label{sec:Projective}

    In this section, we show that the synchronous version of projective splitting in \cite{Combettes2018Asynchronous} is a special case of the four-operator splitting method in Algorithm~\ref{alg:FOS} with $D=E=0$. The algorithm in \cite{Combettes2018Asynchronous} has no step-size restrictions besides being upper and lower bounded. We give an alternative proof to this fact with an appeal to Corollary~\ref{cor:FOS_conv} with $D=E=0$.
    
    We consider a special case of the synchronous version of projective splitting in \cite{Combettes2018Asynchronous} that solves monotone inclusion problems of the form
\begin{align}
  0\in A_n(x)+\sum_{i=1}^{n-1}L_i^*A_i(L_ix)
  \label{eq:PSprob}
\end{align}
that satisfy the following assumption.
\begin{ass}
Assume that $A_i:\mathcal{G}_i\to 2^{\mathcal{G}_i}$ for $i\in\{1,\ldots,n-1\}$ and $A_n:\hilbert\to 2^\hilbert$ are maximally monotone and that $L_i:\hilbert\to\mathcal{G}_i$ for $i\in\{1,\ldots,n-1\}$ are bounded linear operators. Further assume that the solution set to \eqref{eq:PSprob} is nonempty.
  \label{ass:PSprob}
\end{ass}
A synchronous version of the projective splitting algorithm in \cite[Algorithm~3.4]{Combettes2018Asynchronous} reads as follows, where $J_{\gamma A}:=(\id+\gamma A)^{-1}$ is a standard resolvent notation.
\ifx\one\showalgs\begin{algorithm}[H]
  \begin{algorithmic}[1]
    \STATE {\bf Input:} $x_0\in\hilbert$ and $w_{i,0}\in\mathcal{G}_i$ for $i=1,\ldots,n-1$
    \FOR{$k=0,1,\ldots$}
  \STATE $\hat{x}_{k} := J_{\tau_{n,k} A_i}(x_k-\tau_{n,k}\sum_{i=1}^{n-1}L_i^*w_{i,k})$\label{upd:hatxk}
  \STATE $\hat{y}_{k} := (\tau_{n,k}^{-1}x_k-\sum_{i=1}^{n-1}L_i^*w_{i,k})-\tau_{n,k}^{-1}\hat{x}_{k}$\label{upd:hatyk}
    \FOR{$i=1,\ldots,n-1$}
  \STATE $\hat{v}_{i,k} := J_{\tau_{i,k} A_i}(L_ix_k+\tau_{i,k}w_{i,k})$\label{upd:hatvk}
  \STATE $\hat{w}_{i,k} := w_{i,k}+\tau_{i,k}^{-1}L_ix_k-\tau_{i,k}^{-1}\hat{v}_{i,k}$\label{upd:hatwk}
  \ENDFOR
  \STATE $t_k^*:=\hat{y}_{k}+\sum_{i=1}^{n-1}L_i^*\hat{w}_{i,k}$
  \STATE $t_{i,k}:= \hat{v}_{i,k}-L\hat{x}_{k}$
  \STATE $\mu_k:=\frac{\left(\sum_{i=1}^{n-1}\langle t_{i,k},w_{i,k}\rangle-\langle \hat{v}_{i,k},\hat{w}_{i,k}\rangle\right)+\langle t^*,x_k\rangle-\langle \hat{y}_{k},\hat{x}_{k}\rangle}{\sum_{i=1}^{n-1}\|t_{i,k}\|^2+\|t_k^*\|^2}$
  \FOR{$i=1,\ldots,n-1$}
  \STATE $w_{i,k+1} = w_{i,k}-\theta_k\mu_kt_{i,k}$
  \ENDFOR
  \STATE $x_{k+1}:=  x_k-\theta_k\mu_kt_k^*$
		\ENDFOR
	\end{algorithmic}
        \caption{Synchronous Version of Projective Splitting in \cite[Algorithm~3.4]{Combettes2018Asynchronous}}
        \label{alg:PS}
\end{algorithm}\fi

We will reformulate \eqref{eq:PSprob} and show that the four-operator splitting algorithm in Algorithm~\ref{alg:FOS} with $D=E=0$ applied to the reformulation gives Algorithm~\ref{alg:PS}. For the reformulation, we introduce dual variables $w_i\in A_i(L_ix)$ for $i=1,\ldots,n-1$, which is equivalent to that $0\in A_i^{-1}(w_i)-L_ix$. It is straightforward to verify that \eqref{eq:PSprob} is equivalent to the primal-dual formulation
\begin{align}
  0\in \underbrace{\begin{bmatrix}A_{1}^{-1}(w_1)\\\vdots\\ A_{n-1}^{-1}(w_{n-1})\\A_n(x)\end{bmatrix}}_{B}+\underbrace{\begin{bmatrix}0&\ldots& 0 &-L_1\\\vdots&\ddots&\vdots&\vdots\\0&\ldots& 0 &-L_{n-1}\\L_1^*&\ldots& L_{n-1}^* & 0\end{bmatrix}}_{K}\begin{bmatrix}w_1\\\vdots\\w_{n-1}\\x\end{bmatrix}.
  \label{eq:PD-inclusion}
\end{align}
Let $w=(w_1,\ldots,w_{n-1})$ and $p=(w,x)$, to write this as the monotone inclusion problem 
\begin{align}
  0\in Bp+Kp,
\end{align}
where $B$ is maximally monotone since all individual block-operators are maximally monotone and $K$ is linear skew-adjoint, hence maximally monotone. Since $K$ has full domain, also the sum is maximally monotone \cite[Corollary~25.5]{Bauschke2017Convex}. Further, the solution set is nonempty due to Assumption~\ref{ass:PSprob}.

The synchronous projective splitting method stated in Algorithm~\ref{alg:PS} is obtained from Algorithm~\ref{alg:FOS} by letting $D=E=0$, $S=\id$, and 
\begin{align} 
  Q_k := \begin{bmatrix}\tau_{1,k}\id &&&\\&\ddots&&\\&&\tau_{n-1,k}\id&\\&&&\tau_{n,k}^{-1}\id\end{bmatrix},
                                                                                \label{eq:Qdef}
\end{align}
where $\tau_{i,k}>0$ are positive step-sizes. Algorithm~\ref{alg:FOS} then reads as follows with $B$ and $K$ from \eqref{eq:PD-inclusion}, where we have used that $K$ is linear skew-adjoint in the update formula of $\mu_k$.
\ifx\one\showalgs\begin{algorithm}[H]
  \begin{algorithmic}[1]
    \STATE {\bf Input:} $p_0\in\mathcal{G}_1\times\cdots\times\mathcal{G}_{n-1}\times\hilbert$
		\FOR{$k=0,1,\ldots$}
                \STATE $\hat{p}_{k}:=(Q_k+B)^{-1}(Q_k-K)p_{k}$ 
                \STATE $\mu_k:=\frac{\|p_k-\hat{p}_k\|_{Q_k}^2}{\|(Q_k-K)p_k-(Q_k-K)\hat{p}_k\|^2}$
                \STATE $p_{k+1}:=p_k-\theta_k\mu_k((Q_k-K)p_k-(Q_k-K)\hat{p}_k)$
		\ENDFOR
	\end{algorithmic}
        \caption{Synchronous Projective Splitting -- Resolvent Formulation}
        \label{alg:PS_res}
      \end{algorithm}\fi




In Appendix~\ref{app:projective}, we show that Algorithm~\ref{alg:PS_res} is equivalent to Algorithm~\ref{alg:PS}. Convergence of Algorithm~\ref{alg:PS_res}, and hence of synchronous projective splitting in Algorithm~\ref{alg:PS}, follows from the definition of $Q_k$ in \eqref{eq:Qdef} and Corollary~\ref{cor:FOS_conv}. 
    \begin{prop}
      Suppose that Assumption~\ref{ass:PSprob} holds, that $\theta_k\in[\epsilon_\theta,2-\epsilon_\theta]$ for some $\epsilon_\theta\in(0,1)$, and that $\tau_{i,k}\in[\epsilon,\tfrac{1}{\epsilon}]$ for some $\epsilon\in(0,1)$.
      Then Algorithm~\ref{alg:PS_res} generates a sequence $(p_k)_{k\in\naturals}$ such that $p_k\weakto p\in\zero(B+K)$, where $B$ and $K$ are defined in \eqref{eq:PD-inclusion}. Moreover, if $B+K$ is metrically subregular (Definition~\ref{def:metric_subreg}) at all $z\in\zero(B+K)$ for 0 and $\hilbert$ is finite-dimensional, then $x_k\to\bar{x}\in\zero(B+K)$ locally R-linearly and $\dist_S(x_k,\zero(B+K))$ converges locally Q-linearly to zero.
      \label{prop:PSconv}
    \end{prop}
    \begin{pf}
      The result follows from Corollary~\ref{cor:FOS_conv} by verifying the assumptions on $Q_k$. The $Q_k$ in this setting is defined in \eqref{eq:Qdef}. Single-valuedness and maximal monotonicity are obvious. Lipschitz continuity follows from $\tau_{i,k}\in[\epsilon,\tfrac{1}{\epsilon}]$. Let $P=\epsilon\id$, then 1-strong monotonicity w.r.t. $\|\cdot\|_P$ follows from $\tau_{i,k}\geq\epsilon$.
    \end{pf}
    The linear convergence result under metric subregularity appears to be new and uses a weaker assumption than for the linear rates in \cite{Johnstone2019Convergence}.

    \begin{rem}
      Synchronous projective splitting is similar to both long-step FBF (Algorithm~\ref{alg:FOS_gammak} with $E=K=0$) and Chambolle-Pock. Projective splitting is Algorithm~\ref{alg:FOS_gammak} applied to \eqref{eq:PD-inclusion} with $D=E=0$ which gives $M=Q-K$, and projection metric $\|\cdot\|$. This gives the same algorithmic structure as the long-step FBF method, since $D$ and $K$ are handled in the same way in Algorithm~\ref{alg:FOS_gammak} (they only contribute differently to step-size restrictions). To compare to Chambolle-Pock, we use $n=2$ operators in \eqref{eq:PD-inclusion}. 
      Chambolle-Pock can be applied to this problem and is obtained from NOFOB by letting
      \begin{align*}
        M = \underbrace{\begin{bmatrix}\tau_{1}\id&0\\0&\tau_{2}^{-1}\id\end{bmatrix}}_Q-\underbrace{\begin{bmatrix}0&-L_1\\L_1^*&0\end{bmatrix}}_K+\underbrace{\begin{bmatrix}0&0\\2L_1^*&0\end{bmatrix}}_{Z}=\underbrace{\begin{bmatrix}\tau_{1}\id&0\\2L_1^*&\tau_{2}^{-1}\id\end{bmatrix}}_{Q+Z}-\underbrace{\begin{bmatrix}0&-L_1\\L_1^*&0\end{bmatrix}}_K.
      \end{align*}
      This differs from the projective splitting (and the long-step FBF applied to \eqref{eq:PD-inclusion}) kernel in the addition of $Z$. This addition has a few implications.
      \begin{enumerate}[(i)]
      \item The linear kernel $M$ becomes symmetric. This implies that the projection metric kernel $S$ can be chosen equal to $M$. This choice leads to Chambolle-Pock via the FBS method and has the consequence that the second application of $M$ is avoided, see Section~\ref{sec:FBS}. For projective splitting, $M$ is not symmetric. Therefore, the projection kernel $S$ cannot be chosen equal to $M$ and the second application of $M$ is needed to compute the projection. This may be a reason why there is a proof for only this symmetric case in \cite{Chambolle2011AFirstOrder}, although non-symmetric versions are included in the algorithm description. Addition of a projection step would give convergence of the non-symmetric methods. 
      \item The resolvent step in Chambolle-Pock is
      \begin{align*}
        (M+B+K)^{-1}Mp=(Q+Z+B)^{-1}(Q+Z-K)p,
      \end{align*}
      which leads to a Gauss-Seidel serial update between the two blocks in $B$ due to the block-lower triangular structure of $Q+Z$. In projective splitting, we have a similar resolvent update, but without $Z$. Since $Q$ is block-diagonal, the blocks can be updated in parallel.
    \item Finally, the addition of $Z$ puts restrictions on the step-lengths $\tau_1$ and $\tau_2$. The requirement that $M$ should be strongly monotone is, e.g., satisfied by the standard condition for Chambolle-Pock that $\tau_1^{-1}\tau_2\|L_1\|^2<1$ (in our notation). In the projective splitting setting, without the $Z$, there are no restrictions on $\tau_1>0$ and $\tau_2>0$, since $M=Q-K$ is strongly monotone whenever $Q$ is.

\end{enumerate}

    \end{rem}

\fi

    \section{Conclusions}
\label{sec:Conclusions}

We have presented the versatile and conceptually simple algorithm NOFOB (nonlinear forward-backward splitting) for solving monotone inclusion problems. It is based on a novel construction that we call nonlinear forward-backward map. NOFOB is a separate and project method. A nonlinear forward-backward step creates a separating hyperplane between the current point and the solution set. This is followed by a projection onto the hyperplane. We have proven weak convergence to the solution set of the monotone inclusion, and strong and linear convergence under a metric subregularity assumption. We have also presented a novel four-operator splitting method based on NOFOB and shown that many algorithms are special cases of this framework such as forward-backward(-half)-forward splitting, asymmetric forward-backward adjoint splitting with its many special cases, as well as synchonous projective splitting.

\section{Acknowledgments}

We thank the anonomous reviewers for valuable feedback that has greatly improved the manuscript.

\bibliographystyle{siamplain}
\bibliography{references}

\appendix

\section{Defered proofs}
\label{app:def_proofs}
This appendix collects defered proofs.

\subsection{Proof of Proposition~\ref{prop:muk_bounds}}
\label{app:muk_bounds_pf}


We first note that by strong monotonicity of $M_k$, Proposition~\ref{prop:sm_lower} implies that $\|M_kx-M_ky\|_{S^{-1}}=0$ if and only if $x=y$. Therefore, the denominator in the definition of $\mu_{k}^{x,y}$ is positive whenever $x\neq y$.
   
We first prove the lower bound for $x\neq y$. Due to 1-strong monotonicity of $M_k$ w.r.t. $\|\cdot\|_P$ we have for all $x\neq y$
   \begin{align*}
\mu_k^{x,y}&=\frac{\langle M_kx-M_ky,x-y\rangle-\tfrac{\beta}{4}\|x-y\|_P^2}{\|M_kx-M_ky\|_{S^{-1}}^2}
                   \geq(1-\tfrac{\beta}{4})\frac{\|x-y\|_P^2}{\|M_kx-M_ky\|_{S^{-1}}^2}\\
&\geq (1-\tfrac{\beta}{4})\tfrac{\lambda_{\min}(P)}{\lambda_{\max}(S^{-1})}\frac{\|x-y\|^2}{\|M_kx-M_ky\|^2}\geq(1-\tfrac{\beta}{4})\tfrac{\lambda_{\min}(P)}{L_M^2\lambda_{\max}(S^{-1})}>0,
\end{align*}
where the second to last inequality is due to $L_M$ Lipschitz continuity of $M_k$ and since $x\neq y$, and the last holds since $\beta\in[0,4)$ by Assumption~\ref{ass:prob} and $S,P\in\PD$.

For the upper bound, Cauchy-Schwarz inequality and $\beta\geq 0$ imply for all $x\neq y$ that
   \begin{align*}
\mu_k^{x,y}&=\frac{\langle M_kx-M_ky,x-y\rangle-\tfrac{\beta}{4}\|x-y\|_P^2}{\|M_kx-M_ky\|_{S^{-1}}^2}
                   \leq\frac{\|M_kx-M_ky\|_{S^{-1}}\|x-y\|_S}{\|M_kx-M_ky\|_{S^{-1}}^2}\\
     &= \frac{\|x-y\|_S}{\|M_kx-M_ky\|_{S^{-1}}}
       \leq \tfrac{\lambda_{\max}(S)}{\lambda_{\min}(P)},
   \end{align*}
   where 1-strong monotonicity of $M_k$ w.r.t. $\|\cdot\|_P$ and Proposition~\ref{prop:sm_lower} have been used in the second inequality.

Since $\mu_k=\mu_k^{x_k,\hat{x}_k}$ and $x_k\neq\hat{x}_k$ if and only if $x_k\not\in\zero(A+C)$ (Proposition~\ref{prop:fixedFB-zeroAC}), the bounds hold for $\mu_k$ whenever $x_k\not\in\zero(A+C)$.

   \subsection{Proof of Lemma~\ref{lem:M_prop_FOS}}

   \label{app:M_prop_FOS_pf}
   {\it Single-valuedness}: The kernel $M_k=\gamma_k^{-1}\id-D-K$ is single-valued since the three involved operators are. 
   
  {\it Strong monotonicity}: It holds for all $k\in\naturals$ that
\begin{align*}
  \langle M_kx-M_ky,x-y\rangle&=\gamma_k^{-1}\|x-y\|^2-\langle Dx-Dy,x-y\rangle-\langle Kx-Ky,x-y\rangle\\
                              &=\gamma_k^{-1}\|x-y\|^2-\langle Dx-Dy,x-y\rangle\\
                              &\geq\gamma_k^{-1}\|x-y\|^2-\|Dx-Dy\|\|x-y\|\\
                              &\geq\gamma_k^{-1}\|x-y\|^2-L_D\|x-y\|^2\\
                              &\geq\underline{\gamma}^{-1}\|x-y\|^2-L_D\|x-y\|^2\\
                              &=\|x-y\|_P^2,
\end{align*}
where skewness of $K$, Cauchy-Schwarz inequality, $L_D$-Lipschitz continuity of $D$, and the definition of $P$ have been used in the inequalities. This shows that all $M_k$ are 1-strongly monotone w.r.t. $\|\cdot\|_P$. 

{\it Lipschitz continuity}: We have for all $k\in\naturals$:
\begin{align*}
  \|M_kx-M_ky\|&\leq\gamma_k^{-1}\|x-y\|+\|Dx-Dy\|+\|Kx-Ky\|\\&\leq(\gamma_k^{-1}+L_D+\|K\|)\|x-y\|.
\end{align*}
So $M_k$ is $(\gamma_k^{-1}+L_D+\|K\|)$-Lipschitz continuous.

{\it Maximal monotonicity}: Since all $M_k:\hilbert\to\hilbert$ are monotone and (Lipschitz) continuous, they are maximally monotone \cite[Corollary~20.28]{Bauschke2017Convex}.

\subsection{Proof of Proposition~\ref{prop:mu_bound_FOS}}
\label{app:mu_bound_FOS}
Let us define $N:=L_D\id-D$. Since $D$ and $-D$ are $L_D$-Lipschitz continuous, $N$ is $\tfrac{1}{2L_D}$-cocoercive (follows from \cite[Definition~4.10]{Bauschke2017Convex}). Further, $M=N-K+(\gamma^{-1}-L_D)\id$ and
\begin{align}
\nonumber\|&Mx-My\|^2\\
\nonumber    &=\|Nx-(K-(\gamma^{-1}-L_D)\id)x-(Ny-(K-(\gamma^{-1}-L_D)\id)y)\|^2\\
\nonumber             &=\|Nx-Ny\|^2-2\langle Nx-Ny,Kx-Ky-(\gamma^{-1}-L_D)(x-y)\rangle\\
\nonumber             &\quad+\|(K-(\gamma^{-1}-L_D)\id)(x-y)\|^2\\
\nonumber    &=\|Nx-Ny\|^2-2L_D\langle Nx-Ny,x-y\rangle-2\langle Nx-Ny,Kx-Ky-\gamma^{-1}(x-y)\rangle\\
\nonumber  &\quad+\|(K-(\gamma^{-1}-L_D)\id)(x-y)\|^2\\
\nonumber             &\leq-2\langle Nx-Ny,Kx-Ky-\gamma^{-1}(x-y)\rangle+\|(K-(\gamma^{-1}-L_D)\id)(x-y)\|^2\\
\nonumber             &=-2\langle Nx-Ny,Kx-Ky\rangle+2\gamma^{-1}\langle Nx-Ny,x-y\rangle+\|(K-(\gamma^{-1}-L_D)\id)(x-y)\|^2\\
\nonumber             &=2\langle Dx-Dy,Kx-Ky\rangle+2\gamma^{-1}\langle Nx-Ny,x-y\rangle+\|(K-(\gamma^{-1}-L_D)\id)(x-y)\|^2\\
\nonumber    &=2\langle Dx-Dy,Kx-Ky\rangle+2\gamma^{-1}\langle Nx-Ny,x-y\rangle+\|K(x-y)\|^2\\
\nonumber  &\quad+(\gamma^{-1}-L_D)^2\|x-y\|^2\\
\nonumber    &=2\langle Dx-Dy,Kx-Ky\rangle+2\gamma^{-1}\langle Mx-My,x-y\rangle+\|K(x-y)\|^2\\
\nonumber  &\quad+(2\gamma^{-1}(L_D-\gamma^{-1})+(\gamma^{-1}-L_D)^2)\|x-y\|^2\\
\nonumber    &\leq 2\gamma^{-1}\langle Mx-My,x-y\rangle\\
 \nonumber   &\quad+(2L_D\|K\|+\|K\|^2+2\gamma^{-1}(L_D-\gamma^{-1})+(\gamma^{-1}-L_D)^2)\|x-y\|^2,\\\label{eq:FOS_M_ineq_pf}
\end{align}
where cocoercivity of $N$ has been used in the first inequality, skew-adjointness of $K$ in the 5\ordth and 6\ordth equalities, and Cauchy-Schwarz inequality, Lipschitz continuity of $D$, and the definition of $\|K\|$ in the last inequality.

Now, let us show that $\delta$ in \eqref{eq:delta_def} satisfies $\delta\in(0,1)$. We first note that $\bar{\epsilon}$ in \eqref{eq:delta_def} satisfies $\bar{\epsilon}>0$ since $\epsilon\in(0,1)$ and that $\bar{\epsilon}\leq 1$ since
\begin{align*}
  \bar{\epsilon}&=\frac{\epsilon((8-\epsilon)\sqrt{\beta_E^2+16(L_D+\|K\|)^2}+\epsilon\beta_E)}{4(8-\epsilon)}\leq\frac{\epsilon((8-\epsilon)((4-\epsilon)\gamma^{-1}-\beta_E)+\epsilon\beta_E)}{4(8-\epsilon)}\\
  &\leq \frac{\epsilon(8-\epsilon)(4-\epsilon)\gamma^{-1}}{4(8-\epsilon)}= \frac{\epsilon(4-\epsilon)\gamma^{-1}}{4}\leq \epsilon\gamma^{-1}\leq 1,
\end{align*}
where the first inequality comes from rearranging $\gamma\leq\tfrac{4-\epsilon}{\beta_E+\sqrt{\beta_E^2+16(L_D+\|K\|)^2}}$, the second from $\epsilon\in(0,1)$, and the last from $\gamma\geq\epsilon$. Therefore
\begin{align*}
  \delta=\frac{\bar{\epsilon}}{2(\epsilon^{-1}+L_D+\|K\|-\tfrac{\beta_EL_D\epsilon}{4(1-\epsilon)})}\leq \frac{\bar{\epsilon}}{2(\epsilon^{-1}-\tfrac{\beta_EL_D\epsilon}{4(1-\epsilon)})}\leq \frac{\bar{\epsilon}}{\epsilon^{-1}}=\bar{\epsilon}\epsilon\leq\epsilon< 1,
\end{align*}
where the first inequality holds since $L_D+\|K\|\geq 0$ and the second inequality holds since we assume $\tfrac{\epsilon^{-1}}{2}\geq\tfrac{\beta_EL_D\epsilon}{4(1-\epsilon)}$. That $\delta>0$ follows directly from that $\tfrac{\epsilon^{-1}}{2}\geq\tfrac{\beta_EL_D\epsilon}{4(1-\epsilon)}$, which makes the denominator positive.

Let us continue \eqref{eq:FOS_M_ineq_pf} by using $\delta\in(0,1)$ and strong monotonicity of $M$ that makes $\langle Mx-My,x-y\rangle\geq 0$, followed by Cauchy-Schwarz inequality combined with $(\gamma^{-1}+L_D+\|K\|)$-Lipschitz continuity of $M$ (Lemma~\ref{lem:M_prop_FOS}), and thereafter $P=(\gamma^{-1}-L_D)\id$ to conclude
\begin{align*}
  \|Mx-My\|^2&\leq 2\gamma^{-1}\langle Mx-My,x-y\rangle\\
             &\quad+(2L_D\|K\|+\|K\|^2+2\gamma^{-1}(L_D-\gamma^{-1})+(\gamma^{-1}-L_D)^2)\|x-y\|^2\\
  &\leq (2-\delta)\gamma^{-1}\langle Mx-My,x-y\rangle+\delta\gamma^{-1}(\gamma^{-1}+L_D+\|K\|)\|x-y\|^2\\
             &\quad+(2L_D\|K\|+\|K\|^2+2\gamma^{-1}(L_D-\gamma^{-1})+(\gamma^{-1}-L_D)^2)\|x-y\|^2\\
  &=(2-\delta)\gamma^{-1}\langle Mx-My,x-y\rangle+\tfrac{1}{\gamma^{-1}-L_D}\Big(\delta\gamma^{-1}(\gamma^{-1}+L_D+\|K\|)\\
             &\quad+(2L_D\|K\|+\|K\|^2+2\gamma^{-1}(L_D-\gamma^{-1})+(\gamma^{-1}-L_D)^2)\Big)\|x-y\|_P^2.
\end{align*}
Now, if 
\begin{align}
\nonumber  2L_D\|K\|+\|K\|^2+&2\gamma^{-1}(L_D-\gamma^{-1})+(\gamma^{-1}-L_D)^2\\
  &+\delta\gamma^{-1}(\gamma^{-1}+L_D+\|K\|)
  \leq-(2-\delta)\gamma^{-1}(\gamma^{-1}-L_D)\tfrac{\beta_E}{4},\label{eq:param_req}
\end{align}
we conclude that
\begin{align*}
  \|Mx-My\|^2&\leq (2-\delta)\gamma^{-1}(\langle Mx-My,x-y\rangle-\tfrac{\beta_E}{4}\|x-y\|_P^2).
\end{align*}
The left-hand side is 0 if and only if $x=y$ due to strong montonocity of $M$ and Propsosition~\ref{prop:sm_lower}. Therefore, if \eqref{eq:param_req} holds, also the result in \eqref{eq:gamma_over_twodelta} holds for all $x\neq y$ by dividing both sides by $\|Mx-My\|^2> 0$ and $(2-\delta)\gamma^{-1}>0$.

We will finish the proof by proving that \eqref{eq:param_req} holds for all allowed $\gamma$. Since the terms involving $\delta$ in \eqref{eq:param_req} satisfy
\begin{align*}
  \delta\gamma^{-1}((\gamma^{-1}+L_D+\|K\|)-\tfrac{\beta_E(\gamma^{-1}-L_D)}{4})&=\tfrac{\bar{\epsilon}\gamma^{-1}}{2}\tfrac{\gamma^{-1}+L_D+\|K\|-\tfrac{\beta_E(\gamma^{-1}-L_D)}{4}}{\epsilon^{-1}+L_D+\|K\|-\tfrac{\beta_EL_D\epsilon}{4(1-\epsilon)}}\\
                                                                                &\leq\tfrac{\bar{\epsilon}\gamma^{-1}}{2}\tfrac{\epsilon^{-1}+L_D+\|K\|-\tfrac{\beta_E(\tfrac{L_D}{1-\epsilon}-L_D)}{4}}{\epsilon^{-1}+L_D+\|K\|-\tfrac{\beta_EL_D\epsilon}{4(1-\epsilon)}}
                                                                                  \\
                                                                                &=\tfrac{\bar{\epsilon}\gamma^{-1}}{2}\tfrac{\epsilon^{-1}+L_D+\|K\|-\tfrac{\beta_EL_D\epsilon}{4(1-\epsilon)}}{\epsilon^{-1}+L_D+\|K\|-\tfrac{\beta_EL_D\epsilon}{4(1-\epsilon)}}
  =\tfrac{\bar{\epsilon}\gamma^{-1}}{2}
\end{align*}
where the inequality comes from $\gamma\geq\epsilon$ and $\gamma\leq\tfrac{1-\epsilon}{L_D}$, we can verify \eqref{eq:param_req} by verifying
\begin{align}
  2L_D\|K\|+\|K\|^2+&2\gamma^{-1}(L_D-\gamma^{-1})+(\gamma^{-1}-L_D)^2
                      \leq-2\gamma^{-1}(\tfrac{\beta_E}{4}+\tfrac{\bar{\epsilon}}{4}).
                      \label{eq:param_req_simpl}
\end{align}
Multiplying by $2\gamma^2>0$ and rearranging gives
\begin{align*}
  0&\geq \gamma^2(2(2L_D\|K\|+\|K\|^2)+\gamma(4L_D+\beta_E+\bar{\epsilon})-4+2\gamma^2(\gamma^{-2}-2\gamma^{-1}L_D+L_D^2)\\
   &=\gamma^2(2(2L_D\|K\|+\|K\|^2+L_D^2)+\gamma(\beta_E+\bar{\epsilon})-2\\
   &=\gamma^22(L_D+\|K\|)^2+\gamma(\beta_E+\bar{\epsilon})-2,
\end{align*}
which holds if and only if $\gamma$ is between the two roots of the r.h.s. (since the coefficient for $\gamma^2$ is positive), namely
\begin{align*}
  \gamma\in -\frac{\beta_E+\bar{\epsilon}}{4(L_D+\|K\|)^2}+\sqrt{\frac{(\beta_E+\bar{\epsilon})^2}{16(L_D+\|K\|)^4}+\frac{2}{2(L_D+\|K\|)^2}}[-1,1].
\end{align*}
The lower bound is negative, so the inequality holds for all $\gamma>0$ such that
\begin{align}
\nonumber  \gamma &\leq -\frac{\beta_E+\bar{\epsilon}}{4(L_D+\|K\|)^2}+\sqrt{\frac{(\beta_E+\bar{\epsilon})^2}{16(L_D+\|K\|)^4}+\frac{2}{2(L_D+\|K\|)^2}}\\
\nonumber  &= \frac{\sqrt{(\beta_E+\bar{\epsilon})^2+16(L_D+\|K\|)^2}-(\beta_E+\bar{\epsilon})}{4(L_D+\|K\|)^2}\\
         &=\frac{4}{\beta_E+\bar{\epsilon}+\sqrt{(\beta_E+\bar{\epsilon})^2+16(L_D+\|K\|)^2}}.\label{eq:gamma_bound}
\end{align}
The proof is concluded by verifying, with straight-forward but somewhat tedious algebra that is omitted, that the choice of $\bar{\epsilon}$ in \eqref{eq:ep_bar_def} implies that the upper bound for $\gamma$ in \eqref{eq:gamma_bound} that makes \eqref{eq:param_req_simpl} hold satisfies
\begin{align*}
  \frac{4}{\beta_E+\bar{\epsilon}+\sqrt{(\beta_E+\bar{\epsilon})^2+16(L_D+\|K\|)^2}}=\frac{4-\epsilon}{\beta_E+\sqrt{\beta_E^2+16(L_D+\|K\|)^2}}.
\end{align*}
So \eqref{eq:param_req_simpl} and therefore \eqref{eq:param_req} holds for all allowed step-sizes $\gamma$, and the proof is complete.

\section{Projective Splitting Equivalence}
\label{app:projective}

\paragraph{The resolvent step} We start by considering the (nonsymmetric but linear) resolvent step in Algorithm~\ref{alg:PS_res}:
\begin{align}
\nonumber  \hat{p}_k&=(Q_k+B)^{-1}(Q_k-K)p_k\\
\nonumber     &=\left[\begin{array}{l}
       (\tau_{1,k}\id+A_1^{-1})^{-1}(\tau_{1,k}w_{1,k}+L_1x_k)\\
       \qquad\qquad\vdots\\
       (\tau_{n-1,k}\id+A_{n-1}^{-1})^{-1}(\tau_{n-1,k}w_{n-1,k}+L_{n-1}x_k)\\
       (\tau_{n,k}^{-1}\id+A_n)^{-1}(\tau_{n,k}^{-1}x_k-\sum_{i=1}^{n-1}L_i^*w_{i,k})\\
     \end{array}\right]\\
\label{eq:rplus}     &=\left[\begin{array}{l}
       J_{\tau_{1,k}^{-1}A_1^{-1}}(w_{1,k}+\tau_{1,k}^{-1}L_1x_k)\\
       \qquad\qquad\vdots\\
       J_{\tau_{n-1,k}^{-1}A_{n-1}^{-1}}(w_{n-1,k}+\tau_{n-1,k}^{-1}L_{n-1}x_k)\\
       J_{\tau_{n,k} A_n}(x_k-\tau_{n,k}\sum_{i=1}^{n-1}L_i^*w_{i,k})\\
     \end{array}\right].
\end{align}
To arrive at the updates in Algorithm~\ref{alg:PS}, we make use of Moreau's identity \cite[Proposition~23.20]{Bauschke2017Convex} 
\begin{align}
  J_{\tau A}(z)=z-\tau J_{\tau^{-1}A^{-1}}(\tau^{-1}z) 
  \label{eq:Moreau}
\end{align}
and store some additional variables. For the $n-1$ first blocks in \eqref{eq:rplus}, we let in \eqref{eq:Moreau}; $A=A_i$, $\tau=\tau_{i,k}$, and resolvent input $z=\tau_{i,k}(w_{i,k}+\tau_{i,k}^{-1}L_ix_k)$. Then $\hat{w}_{i,k}=J_{\tau_{i,k}^{-1}A_{i}^{-1}}(w_{i,k}+\tau_{i,k}^{-1}L_{i}x_k)$ in \eqref{eq:rplus} can be computed as
\begin{align}
\label{eq:hatvk}\hat{v}_{i,k} &:= J_{\tau_{i,k} A_i}(L_ix_k+\tau_{i,k}w_{i,k}),\\
\label{eq:hatwk}\hat{w}_{i,k} &:= w_{i,k}+\tau_{i,k}^{-1}L_ix_k-\tau_{i,k}^{-1}\hat{v}_{i,k}.
\end{align}
These are the updates on Lines~\ref{upd:hatvk}~and~\ref{upd:hatwk} in Algorithm~\ref{alg:PS}. It is straightforward to show that $(\hat{v}_{i,k},\hat{w}_{i,k})\in\gph(A_i)$, which is crucial for the convergence analysis in \cite{Combettes2018Asynchronous}. We also store an additional point $y_k=J_{\tau_{n,k}^{-1}A_n^{-1}}(\tau_{n,k}^{-1}x_k-\sum_{i=1}^{n-1}L_i^*w_{i,k})$ (which is not needed but used to compare to the update in Algorithm~\ref{alg:PS}) to compute a point $(\hat{x}_k,\hat{y}_k)\in\gph(A_n)$ where $\hat{x}_k=J_{\tau_{n,k} A_n}(x_k-\tau_{n,k}\sum_{i=1}^{n-1}L_i^*w_{i,k})$ is the $n$\ordth block update in \eqref{eq:rplus}. We let in \eqref{eq:Moreau}; $A=A_n$, $\tau=\tau_{n,k}$, and resolvent input $z=x_k-\tau_{n,k}\sum_{i=1}^{n-1}L_i^*w_{i,k}$ to get
\begin{align}
\label{eq:hatxk}\hat{x}_{k} &:= J_{\tau_{n,k} A_n}(x_k-\tau_{n,k}\sum_{i=1}^{n-1}L_i^*w_{i,k}),\\
\label{eq:hatyk}\hat{y}_{k} &:= (\tau_{n,k}^{-1}x_k-\sum_{i=1}^{n-1}L_i^*w_{i,k})-\tau_{n,k}^{-1}\hat{x}_k.
\end{align}
These are the updates on Lines~\ref{upd:hatxk}~and~\ref{upd:hatyk} in Algorithm~\ref{alg:PS}. We have shown that the resolvent update in Algorithm~\ref{alg:PS_res} exactly corresponds to Lines~\ref{upd:hatxk} to \ref{upd:hatwk} in Algorithm~\ref{alg:PS}.


\paragraph{The Projection}

Let us derive an expression for the $\mu_k$ update in Algorithm~\ref{alg:PS_res}. We let $t_k^*:=\hat{y}_k+\sum_{i=1}^{n-1}L_i^*\hat{w}_{i,k}$ and $t_{i,k}:= \hat{v}_{i,k}-L\hat{x}_k$. The numerator satisfies
\begin{align*}
  \|&p_k-\hat{p}_k\|_{Q_k}^2\\
  &=\sum_{i=1}^{n-1}\tau_{i,k}\|w_{i,k}-\hat{w}_{i,k}\|^2+\tau_{n,k}^{-1}\|x_{k}-\hat{x}_k\|^2\\
                           &=\left(\sum_{i=1}^{n-1}\langle \hat{v}_{i,k}-L_ix_k,w_{i,k}-\hat{w}_{i,k}\rangle\right)+\langle \hat{y}_k+\sum_{i=1}^{n-1}L_i^*w_{i,k},x_k-\hat{x}_k\rangle\\
                           &=\left(\sum_{i=1}^{n-1}\langle t_{i,k}-L_i(x_k-\hat{x}_k),w_{i,k}-\hat{w}_{i,k}\rangle\right)+\langle t_k^*+\sum_{i=1}^{n-1}L_i^*(w_{i,k}-\hat{w}_{i,k}),x_k-\hat{x}_k\rangle\\
                           &=\left(\sum_{i=1}^{n-1}\langle t_{i,k},w_{i,k}-\hat{w}_{i,k}\rangle\right)+\langle t_k^*,x_k-\hat{x}_k\rangle\\
                           &=\left(\sum_{i=1}^{n-1}\langle t_{i,k},w_{i,k}\rangle-\langle \hat{v}_{i,k}-L\hat{x}_k,\hat{w}_{i,k}\rangle\right)+\langle t_k^*,x_k\rangle-\langle \hat{y}_k+\sum_{i=1}^{n-1}L_i^*\hat{w}_{i,k},\hat{x}_k\rangle\\
                           &=\left(\sum_{i=1}^{n-1}\langle t_{i,k},w_{i,k}\rangle-\langle\hat{v}_{i,k},\hat{w}_{i,k}\rangle\right)+\langle t_k^*,x_k\rangle-\langle \hat{y}_k,\hat{x}_k\rangle,
\end{align*}
which is exactly the numerator of $\mu_k$ in Algorithm~\ref{alg:PS}.


The denominator of $\mu_k$ in Algorithm~\ref{alg:PS_res} satisfies
\begin{align*}
  \|(Q_k-K)p_k-&\,(Q_k-K)\hat{p}_k\|^2\\
             &=\left(\sum_{i=1}^{n-1}\|\tau_{i,k}w_{i,k}-L_ix_k-(\tau_{i,k}\hat{w}_{i,k}-L_i\hat{x}_k)\|^2\right)\\
             &\quad+\|\tau_{n,k}^{-1}x_k+\sum_{i=1}^{n-1}L_i^*w_{i,k}-(\tau_{n,k}^{-1}\hat{x}_k+\sum_{i=1}^{n-1}L_i^*\hat{w}_{i,k})\|^2\\
             &=\left(\sum_{i=1}^{n-1}\|\hat{v}_{i,k}+L_i\hat{x}_k\|^2\right)+\|\hat{y}_k-\sum_{i=1}^{n-1}L_i^*\hat{w}_{i,k}\|^2\\
  &=\left(\sum_{i=1}^{n-1}\|t_{i,k}\|^2\right)+\|t_k^*\|^2,
\end{align*}
which coincides with the denominator of $\mu_k$ in Algorithm~\ref{alg:PS}. Hence the $\mu_k$ in Algorithm~\ref{alg:PS} and Algorithm~\ref{alg:PS_res} are the same. Using the same equalities, it follows that
\begin{align*}
(Q_k-K)p_k-(Q_k-K)\hat{p}_k=(t_{1,k},\ldots,t_{n-1,k},t_k^*).
\end{align*}
The algorithm update in Algorithm~\ref{alg:PS_res} therefore becomes
\begin{align*}
  p_{k+1} = p_k-\theta_k\mu_k(t_1,\ldots,t_{n-1},t^*),
\end{align*}
which is exactly the update in Algorithm~\ref{alg:PS}. Hence, the algorithms are equivalent.

\end{document}